\documentclass[12pt,a4paper]{amsart}
\usepackage[top=30mm, bottom=30mm, left=25mm, right=25mm]{geometry}
\usepackage{graphicx}
\usepackage{amsthm,amscd,amsmath,amsthm,amssymb,mathrsfs,amsfonts}
\usepackage{appendix}
\usepackage{tikz}
\usepackage{verbatim}
\usepackage{url}
\usepackage[all]{xy}
\usepackage{xcolor}
\usepackage{color}
\usepackage[colorlinks=true,citecolor=blue]{hyperref}
\usepackage{multicol}
\theoremstyle{plain}

\newtheorem{thm}{Theorem}[section]
\newtheorem{lem}[thm]{Lemma}

\theoremstyle{definition}
\newtheorem{defn}[thm]{Definition}

\newtheorem{rem}[thm]{Remark}

\begin{document}
	
\title
[Topological multiple  recurrence for generalized polynomials] 
{Topological multiple  recurrence of weakly mixing minimal systems for generalized polynomials  }

\author[R.~Zhang]{Ruifeng Zhang}
\address[R.~Zhang]{School of Mathematics, Hefei University of Technology, Hefei, Anhui,230009,P.R. China}
\email{rfzhang@hfut.edu.cn}

\author[J.~Zhao]{Jianjie Zhao}
\address[J.~Zhao]{Wu Wen-Tsun Key Laboratory of Mathematics, USTC, Chinese Academy of Sciences and
	School of Mathematics, University of Science and Technology of China,
	Hefei, Anhui, 230026, P.R. China}
\email{zjianjie@mail.ustc.edu.cn}

\date{\today}
\keywords{General polynomials, weakly mixing minimal systems}
\subjclass[2010]{37B20, 37B05}
\maketitle

\maketitle	

\begin{abstract}
	Let $(X, T)$ be a weakly mixing minimal system, $p_1, \cdots, p_d$ be integer-valued generalized polynomials and $(p_1,p_2,\cdots,p_d)$ be non-degenerate.
Then there exists a residual subset $X_0$ of $X$ such that for all $x\in X_0$
$$\{ (T^{p_1(n)}x, \cdots, T^{p_d(n)}x): n\in \mathbb{Z}\}$$
is dense in $X^d$.      
\end{abstract}

\section{Intruduction}

By a topological dynamical system $(X,T)$, we mean a compact metric space $X$  together with  a homeomorphism $T$ from $X$ to itself. By a measure preserving system we mean a quadruple $(X,\mathcal{B},\mu,T)$, where $(X,\mathcal{B},\mu)$ is a Lebesgue space and $T$ and $T^{-1}$ are measure preserving transformations.
In this paper, we study the topological multiple recurrence of  weakly mixing minimal systems.

For a measure preserving system, Furstenberg \cite{Fur} proved the multiple recurrence
theorem, and gave a new proof of Szemer$\acute{e}$di's theorem. Later, Glasner \cite{Glasner} 
considered the counterpart of \cite{Fur} in topological dynamics and proved that: for a weakly mixing minimal system $(X, T)$
and a positive integer $d$,
there is a dense $G_{\delta}$ subset $X_0$ of $X$ such that for each $x\in X_0$,
$\{(T^nx, \cdots, T^{dn}x): n\in \mathbb{Z}\}$ is dense in $X^d$. Note that a different proof of this result can  also be found in \cite{KO,TK}

For a weakly mixing measure preserving system, Bergelson \cite{Be2} proved the following result:
let $(X,\mathcal{B},\mu,T)$ be a weakly mixing system, let $k \in \mathbb{N} $ and let $p_i(n)$ be integer-valued polynomials such that no $p_i$ and no $p_i-p_j$ is constant, $1 \le i\neq j \le k$. Then for any $f_1,f_2,\dots,f_k \in L^{\infty}(X)$,
$$\lim_{N-M \rightarrow \infty}||\frac{1}{N-M}\sum_{n=M}^{N-1}T^{p_1(n)}f_1T^{p_2(n)}f_2\dots T^{p_k(n)}f_k-\prod_{i=1}^{k}\int f d\mu||_{L^2}=0.$$
Note that this is a special case of a polynomial extension of Szemer$\acute{e}$di's theorem obtained in \cite{BL1}.

In the topological side, Huang, Shao and Ye \cite{HSY19} considered the correspondence result of \cite{BL1}, 
and they proved the following result:
let $(X, T)$ be a weakly mixing minimal system and $p_1, \cdots, p_d$ be distinct polynomials with $p_i(0)=0, i=1, \cdots, d$,
then there is a dense $G_{\delta}$ subset $X_0$ of $X$ such that
for each $x\in X_0$,
$$\{ (T^{p_1(n)}x, \cdots, T^{p_d(n)}x): n\in \mathbb{Z}\}$$
is dense in $X^d$.

The multiple recurrence of a weakly mixing measure preserving system for generalized 
polynomials was studied by Bergelson and McCutcheon \cite{BeMc} (for more details concerning generalized polynomials, see \cite{BL2}).
In this paper, we consider the problem in topological side.
As the generalized polynomials are much more complicated than the polynomials, for instance
$\left\lceil 2\pi n -\left\lceil 2\pi n\right\rceil \right\rceil$ can only take values $0$ and $1$, clearly we should preclude such kind of ``bad" generalized polynomials.
So we introduced the notion of $(p_1,p_2,\cdots,p_d)$ be non-degenerate (see Definition \ref{def-nondegenerate}).
The main result of this paper is the following theorem.

\begin{thm} \label{thm general}
	Let $(X, T)$ be a weakly mixing minimal system, $p_1, \cdots, p_d$ be integer-valued generalized polynomials and $(p_1,p_2,\cdots,p_d)$ be non-degenerate.
	Then there is a dense $G_{\delta}$ subset $X_0$ of $X$ such that for all $x\in X_0$,
	$$\{(T^{p_1(n)}x, \cdots, T^{p_d(n)}x): n\in \mathbb{Z}\}$$
	is dense in $X^d$.
	
	Moreover, for any non-empty open subsets $U, V_1, \cdots, V_d$ of $X$, for any $\varepsilon>0$, for any $s, t \in \mathbb{N}$ and $g_1, \cdots, g_t\in \widehat{SGP_s}$, let
	$$C=C(\varepsilon, g_1, \cdots, g_t ):=\bigcap_{k=1}^{t}\{n \in \mathbb{Z}:\{g_k(n)\}\in (-\varepsilon,\varepsilon) \},$$
	$$N=\{n\in \mathbb{Z}: U\cap T^{-p_1(n)}V_1 \cap \cdots \cap T^{-p_d(n)} V_d \neq \emptyset\}.$$
	Then $N \cap C $ is syndetic, where $\widehat{SGP_s}$ and $\{g_k(n)\}$ are defined in Section 2.
\end{thm}

\medskip

The key ingredient in the proof of the main result is to view the integer-valued generalized polynomials, in some sense, as the ordinary polynomials, and thus we can use the method in \cite{HSY19}.
Roughly speaking, the difficulty is in calculating $p(n+m)-p(m)-p(n)$. For instance,  generally $\left \lceil a(n+m)^2 \right \rceil$ is not equal to $\left \lceil an^2 \right \rceil+\left \lceil 2amn \right \rceil+\left \lceil am^2 \right \rceil$, while  $a(n+m)^2=an^2+2anm+am^2$.
To overcome this, we need to restrict $n$ in  some set $C$ where the fractional part $\{an^2\}$ and $\{2amn\}$ are small enough such that for any $n \in C$,
$\left \lceil a(n+m)^2 \right \rceil=\left \lceil an^2 \right \rceil+\left \lceil 2amn \right \rceil+\left \lceil am^2 \right \rceil$.    
Roughly speaking,  we will restrict integer-valued generalized polynomials to a Nil Bohr$_0$-set rather than $\mathbb{Z}$.

The paper is organized as follows.
In Section 2, we  introduce some notions and some properties that will be needed in the proof.
In Section 3, we prove Theorem \ref{thm general} for  integer-valued generalized polynomilals of degree $1$. In the final section, we recall the PET-induction and show the proof of Theorem \ref{thm general}.

\medskip
\noindent {\bf Acknowledgments.}
The authors would like to thank Professor X. Ye  for  help discussions.
The first author were supported by NNSF of China(11871188,11671094), the second author were supported by NNSF of China (11431012).

\section{Preliminary}

\subsection{Some important subsets of integers and Furstenberg families}

In this paper, the set of all integers and positive integers are denoted by $\mathbb{Z}$
and $\mathbb{N}$ respectively.

A subset $S$ of $\mathbb{Z}$ is \emph{syndetic} if it has a bounded gap, i.e.
there is $L\in \mathbb{N}$ such that $\{n, n+1, \cdots, n+L\} \cap S \neq \emptyset$ for every $n\in \mathbb{Z}$.
$S$ is \emph{thick} if it contains arbitrarily long runs of integers, i.e. for any $L\in \mathbb{N}$, there is $a_L\in \mathbb{Z}$ such that $\{a_L, a_L+1, \cdots, a_L+L\}\subset S$.
$S$ is \emph{thickly syndetic} if for every $L\in \mathbb{N}$, there exists a syndetic set $B_L\subset \mathbb{Z}$ such that
$B_L+\{0, 1, \cdots, L\} \subset A$, where $B_L+\{0, 1, \cdots, L\}=\cup_{b\in B_L}\{b, b+1, \cdots, b+L\}$.

The family of all syndetic sets, thick sets and thickly syndetic sets are denoted by $\mathcal{F}_s$, $\mathcal{F}_t$ and $\mathcal{F}_{ts}$ respectively.

Let $\mathcal{P}$ denote the collection of all subsets of $\mathbb{Z}$.
A subset $\mathcal{F}$ of $\mathcal{P}$ is called a \emph{Furstenberg family} (or just a \emph{family}), if it is hereditary upward, i.e.,
$$F_1 \subset F_2 \ \  \text{and} \  \  F_1 \in \mathcal{F} \ \ \text{imply} \ \ F_2\in \mathcal{F}.$$
A family $\mathcal{F}$ is called \emph{proper} if it is a non-empty proper subset of $\mathcal{P}$, i.e. it is neither empty nor all of $\mathcal{P}$.
Any non-empty collection $\mathcal{A}$ of subsets of $\mathbb{Z}$ naturally generates a family
$$\mathcal{F}(A)=\{ F\subset \mathbb{Z}: A\subset F \ \ \text{for some} \ \ A\in \mathcal{A} \}.$$
A proper family $\mathcal{F}$ is called a \emph{filter} if $F_1, F_2 \in \mathcal{F}$ implies
$F_1 \cap F_2 \in \mathcal{F}$.

Note that the set of all thickly syndetic sets is a filter, i.e. the intersection of any finite thickly syndetic sets is still a thickly syndetic set.

\subsection{Topological dynamics}

Let $(X, T)$ be a dynamical system.
For $x\in X$, we denote the orbit of $x$ by
$orb(x, T)=\{T^{n}x : n\in \mathbb{Z}\}$.	A point $x\in X$ is called a \emph{transitive point} if the orbit of $x$ is dense in $X$,
i.e., $\overline{orb(x, T)}=X$.
A dynamical system $(X, T)$ is called \emph{minimal} if every point $x\in X$ is a transitive point.

Let $U, V \subset X$ be two non-empty open sets,
the \emph{hitting time set} of $U$ and $V$ is denoted by
$$N(U, V)=\{n\in \mathbb{Z}: U \cap T^{-n}V \neq \emptyset\}.$$

We say that $(X, T)$ is \emph{(topologically) transitive} if for any non-empty open sets
$U, V \subset X$, the hitting time $N(U, V)$ is non-empty;
\emph{weakly mixing} if the product system $(X\times X, T\times T)$ is transitive.

We say that $(X, T)$ is \emph{thickly syndetic transitive} if for any non-empty open sets
$U, V \subset X$, the hitting time $N(U, V)$ is thickly syndetic.
Let $p_i: \mathbb{Z} \rightarrow \mathbb{Z},i=1,2,\cdots,k$, we say that $(X,T)$ is $\{p_1,p_2,\cdots,p_k\}$-thickly-syndetic transitive if
for any non-empty open sets $U_i, V_i \subset X,i=1,2,\cdots, k$,
$$N(\{p_1,p_2,\cdots, p_k\},U_1\times U_2 \times \cdots \times U_k,V_1 \times V_2 \times \cdots \times V_k ):=\bigcap_{i=1}^{k} N(p_i, U_i,V_i)$$		
is thickly syndetic, where 	$N(p_i, U_i,V_i):=\{n\in \mathbb{Z}: U_i \cap T^{-p_i(n)}V_i \neq \emptyset\}$, $i=1,2,\cdots,k$.

The following Lemma is the analogue	of Lemma $2.6$ in \cite{HSY19}.
\begin{lem}\label{lem-deg-1}
	Let $(X, T)$ be a dynamical system and
	$p_1, \cdots, p_d: \mathbb{Z} \rightarrow \mathbb{Z}$ such that $(X, T)$ is $\{p_1(n), \cdots, p_d(n)\}$-thickly-syndetic transitive.
	Let $C$ be a syndetic set.
	Then for any non-empty open sets $V_1, \cdots, V_d$ of $X$
	and any subsequence $\{ r(n)\}_{n=0}^{\infty}$ of natural numbers,
	there is a sequence of integers $\{k_n\}_{n=0}^{\infty}\subset C$ such that
	$|k_0|> r(0), |k_n|>|k_{n-1}|+r(|k_{n-1}|)$ for all $n \ge 1$,
	and for each $i\in \{1, 2, \cdots, d\}$,
	there is a descending sequence $\{V_{i}^{(n)}\}_{n=0}^{\infty} $ of non-empty open subsets of $V_i$ such that for each $n \ge 0$ one has that
	$$T^{ p_i(k_j)} T^{-j}V_{i}^{(n)} \subset V_i, \ \text{for all \ } 0 \le j \le n. $$
\end{lem}

\begin{proof}
	Let $V_1, \cdots, V_d$ be non-empty open subsets of $X$.
	Then $\bigcap_{i=1}^{d}N(p_i,V_i, V_i )$ is  thickly syndetic.
	Since $C$ is syndetic,
	thus $\bigcap_{i=1}^{d}N(p_i,V_i, V_i ) \cap C$ is syndetic.
	Choose $k_0 \in \bigcap_{i=1}^{d}N(p_i, V_i, V_i ) \cap C$
	such that $|k_0|>r(0)$. It implies
	$ T^{-p_i(k_0)}V_i \cap V_i \neq \emptyset$
	for all $i=1, \cdots, d$.
	Put $V_i^{(0)}=T^{-p_i(k_0)}V_i \cap V_i$ for all  $i=1, \cdots, d$ to complete the base step.
	
	Now assume that for $n\ge 1$ we have found numbers
	$k_0, k_1, \cdots, k_{n-1} \in C$ and for each
	$i=1, \cdots, d$, we have non-empty open subsets
	$V_i \supseteq V_i^{(0)} \supseteq V_i^{(1)} \cdots \supseteq V_i^{(n-1)}$
	such that $|k_0|>r(0)$, and for each $m=1, \cdots, n-1$ one has $|k_m|>|k_{m-1}|+r(|k_{m-1}|)$ and
	$$T^{  p_i(k_j)} T^{-j}V_{i}^{(m)} \subset V_i, \ \text{for \ all \ } 0 \le j \le m.$$
	For $i=1, \cdots, d$, let $U_i=T^{-n}(V_i^{n-1})$.
	Since  $(X, T)$ is $\{p_1(n), \cdots, p_d(n)\}$-thickly-syndetic transitive,
	$$\bigcap_{i=1}^{d}N(p_i, U_i, V_i)=\{n\in \mathbb{Z}: U_i \cap T^{-p_i(n)}V_i\neq \emptyset\} $$
	is thickly syndetic.
	Hence $C \cap (\bigcap_{i=1}^{d}N(p_i,U_i,V_i))$ is syndetic.
	Then there exists $k_n \in C \cap (\bigcap_{i=1}^{d}N(p_i,U_i,V_i))$ such that $|k_n|>|k_{n-1}|+r(|k_{n-1}|)$. It implies
	$$T^{-p_i(k_n) }V_i \cap U_i \neq \emptyset \ \ \text{for all}\ i=1, \cdots, d.$$

	Then for $i=1, \cdots, d$,
	$$T^{p_{i}(k_n)}U_i \cap V_i= T^{p_{i}(k_n)}T^{-n}(V_i^{n-1}) \cap V_i \neq \emptyset.$$
	
	Let
	$$V_i^{(n)}=V_i^{(n-1)}\cap (T^{p_{i}(k_n)}T^{-n})^{-1}V_i.$$
	Then $V_i^{(n)} \subset V_i^{(n-1)}$ is a non-empty open set and
	$$T^{p_{i}(k_n)}T^{-n}V_i^{(n)} \subset V_i. $$
	Since $ V_i^{(n)} \subset V_i^{(n-1)}$,
	we have
	$$ T^{p_{i}(k_j)} T^{-j}V_{i}^{(n)} \subset V_i, \  \text{for all} \ 0 \le j \le n.$$
	Hence we finish our induction. The proof is completed.
\end{proof}	

\medskip

The following Lemma is the analogue of Propostion $1$ in \cite{TK}.	
\begin{lem} \label{lem-equivalent}
	Let $(X, T)$ be a dynamical system and $d\in \mathbb{N}$.
	For any functions $p_1, \cdots, p_d$ from $\mathbb{Z}$ to $\mathbb{Z}$.
	Then the following are equivalent:
	\begin{enumerate}
		\item  If $U, V_1, \cdots, V_d\subset X$ are non-empty open sets,
		then there exists $n\in \mathbb{Z}$, such that
		$$U\cap T^{-p_1(n)}V_1\cap \cdots \cap T^{-p_d(n)}V_d \neq \emptyset.$$
		\item There exists a dense $G_{\delta}$ subset $Y\subset X$ such that
		for every $x\in Y$,
		$$\{(T^{p_1(n)}x, T^{p_2(n)}x, \cdots, T^{p_d(n)}x): n\in \mathbb{Z}\}$$
		is dense in $X^d$.
	\end{enumerate}
\end{lem}	

\begin{proof}
	The proof is similar to the proof in  \cite{TK}. For completeness, we include a proof.
	
	$1 \Rightarrow 2$:
	Consider a countable base of open balls
	$\{B_k: k\in \mathbb{N}\}$ of $X$.
	Put
	$$Y=\bigcap_{(k_1, \cdots, k_d)\in \mathbb{N}^d} \bigcup_{n\in \mathbb{Z}} \bigcap_{i=1}^{d}T^{-p_i(n)}B_{k_i}.$$
	The set $\cup_{n\in \mathbb{Z}} \cap_{i=1 }^{d}T^{-p_i(n)}B_{k_i}$ is open, and is dense by 1.
	Thus by the Baire category theorem, $Y$ is a dense $G_{\delta}$ subset of $X$.
	By construction, for every $x\in Y$, $$\{ (T^{p_1(n)}x, T^{p_2(n)}x, \cdots, T^{p_d(n)}x): n\in \mathbb{Z}\}$$ is dense in $X^d$.
	
	$2 \Rightarrow 1$: Choose $x\in Y\cap U$ and $n\in \mathbb{Z}$ such that
	$$(T^{p_1(n)}x, T^{p_2(n)}x, \cdots, T^{p_d(n)}x) \in V_1\times \cdots \times V _d, $$
	then $x\in U\cap T^{-p_1(n)}V_1\cap \cdots \cap T^{-p_d(n)}V_d $.
\end{proof}

\subsection{Generalized polynomials}
For a real number $a$, let
$\left\| {a} \right\|= \inf\{|a-n| : n\in \mathbb{Z}\}$
and
$ \left\lceil {a} \right\rceil=\min \{m\in \mathbb{Z}: |a-m|=\left\| {a} \right\|\} $.
We denote $[a]$ the greatest integer not exceeding $a$, then $\left\lceil a \right\rceil=[a+\frac{1}{2}]$.
We put $\{a\}=a-\left\lceil {a} \right\rceil$, and $\{a\}\in (-\frac{1}{2},\frac{1}{2}]$.

In \cite{HSY16}, Huang, Shao and Ye introduced the notions of $GP_d$ and $\mathcal{F}_{GP_d}$.

\begin{defn}
	Let $d\in \mathbb{N},$  the \emph{generalized polynomials} of degree $\leq d$ (denoted by $GP_d$) is defined as follows.   	For $d=1$, $GP_1$ is the smallest collection of functions from $\mathbb{Z}$ to $\mathbb{R}$
	containing $\{h_a: a\in \mathbb{R} \} $ with $h_a(n)=an$ for each $n\in \mathbb{Z}$, which is closed under taking $\left\lceil {\cdot } \right\rceil$, multiplying by a constant and finite sums.
	
	Assume that $GP_i$ is defined for $i<d$. Then $GP_d$ is the smallest collection of functions from
	$\mathbb{Z}$ to $\mathbb{R}$ containing $GP_i$ with $i<d$, functions of the forms
	$$a_0n^{p_0}\left\lceil {f_1(n) } \right\rceil \cdots \left\lceil {f_k(n)} \right\rceil $$
	(with $a_0\in \mathbb{R}, p_0 \ge 0, k\ge 0, f_l\in GP_{p_l}$ and $\sum_{l=0}^{k}p_l=d$),
	which is closed under taking $\left\lceil {\cdot } \right\rceil$, multiplying by a constant and finite sums.
	Let $GP=\bigcup_{i=1}^{\infty}GP_i$. Note that if $p\in GP$, then $p(0)=0$.
\end{defn}
\begin{defn}
	Let $\mathcal{F}_{GP_d}$ be the family generated by the sets of forms
	$$\bigcap_{i=1}^{k}\{n\in \mathbb{Z}: p_i(n) \ ( \text{mod}\ \mathbb{Z})\ \in (-\varepsilon_i, \varepsilon_i)\},$$
	where $k\in \mathbb{N}$, $p_i\in GP_d$, and $\varepsilon_i>0, 1\le i\le k$.
	Note that $ p_i(n) \ ( \text{mod} \ \mathbb{Z}) \in (-\varepsilon_i, \varepsilon_i)$
	if and only if $\{p_i(n)\} \in (-\varepsilon_i, \varepsilon_i)$.

\end{defn}
\begin{rem}\label{rem-filer}
	$\mathcal{F}_{GP_d}$ is a filter.
\end{rem}

A subset $A\subset \mathbb{Z}$ is a \emph{Nil$_d$ Bohr$_0$-set} if there exist a $d$-step nilsystem $(X, T)$, $x_0\in X$ and an open set $U\subset X$ containing $x_0$ such that
$N(x_0, U):= \{n\in \mathbb{Z}: T^nx_0 \in U\}$ is contained in $A$.
Denote by $\mathcal{F}_{d, 0}$ the family consisting of all Nil$_d$ Bohr$_0$-sets.  
A subset $A\subset \mathbb{Z}$ is called Nil Bohr$_0$-set if $A\in \mathcal{F}_{d, 0}$ for some $d\in \mathbb{N}$.    	
In \cite{HSY16}, the authors proved the following theorem.
\begin{thm} [Theorem B in \cite{HSY16}] \label{Bohr}
	Let $d\in \mathbb{N}$. Then $ \mathcal{F}_{d, 0}=\mathcal{F}_{GP_d}$.
\end{thm}
\begin{rem}
	Since a nilsystem is distal,  every Nil$_d$ Bohr$_0$-set is syndetic. Together with  Remark \ref{rem-filer} we know $\mathcal{F}_{GP_d}$ is a filter and any $A \in \mathcal{F}_{GP_d}$ is a syndetic set.
\end{rem}

Now we introduce the notion of integer-valued generalized polynomials.
\begin{defn}
	For $d\in \mathbb{N}$, {\it the integer-valued generalized polynomials} of degree $\le d$ is defined by
	$$\widetilde{GP_d}=\{\left\lceil p(n) \right\rceil: p(n)\in GP_d\},$$
	and {\it the integer-valued generalized polynomials}  is then defined by
	$$\mathcal{G}=\bigcup_{i=1}^{\infty}\widetilde{GP_i}.$$
\end{defn}

For $p(n) \in \mathcal{G}$,
the least $d\in \mathbb{N}$ such that $p\in \widetilde{GP_d}$ is defined by the \emph{degree} of $p$, denoted by
$deg(p)$.

\medskip

Since the integer-valued generalized polynomials are very complicated, we will also
specify a subclass of the integer-valued generalized polynomials, i.e. {\it the special integer-valued generalized polynomials} (denoted by $\widetilde{SGP}$),
which will be used in the proof of our main results. 

We  need to recall the defintion of $L(a_1,a_2,\dots,a_l)$ in Defintion 4.2 of  \cite{HSY16}. For $a \in \mathbb{R}$, we define $L(a)=a$. For $a_1,a_2 \in \mathbb{R}$, we define $L(a_1,a_2)=a_1\left\lceil L(a_2)\right\rceil$. Inductively, for $a_1,a_2, \dots, a_l \in \mathbb{R}(l \ge 2)$ we define
$$L(a_1,a_2,\dots,a_l)=a_1\left\lceil L(a_2,\dots,a_l)\right\rceil.$$ 

Before introducing the definition of $\widetilde{SGP}$,
we need to introduce the notion of  the simple generalized polynomials.

\begin{defn} 	
	For $d \in \mathbb{N}$, the {\it simple generalized polynomials} of degree $ \le d$ \  (denoted by $\widehat{SGP_d}$\ ) is defined as follows.
	$\widehat{SGP}_d$ is the smallest collection  of generalized polynomials of the forms
	$$\prod_{i=1}^{k}L(a_{i,1}n^{j_{i,1}},\cdots, a_{i,l_{i}}n^{j_{i,l_i}}),$$
	where $k \ge 1$, $1 \le l_i \le d$,  $a_{i,1},a_{i,2},\dots,a_{i,l_{i}} \in \mathbb{R}$, $j_{i,1},j_{i,2},\dots, j_{i,l_i} \ge 0$ and $\sum_{i=1}^{k}\sum_{t=1}^{l_i}j_{i,t} \le d$.
\end{defn}


With the help of the above definition, we can intoduce the notion of  special integer-valued generalized polynomials.
\begin{defn} 	For $d \in \mathbb{N}$, the {\it special integer-valued generalized polynomials} of degree $ \le d$ (denoted by $\widetilde{SGP_d}$) is defined as follows.
	$$\widetilde{SGP_d}=\{\sum_{i=1}^k c_i\left\lceil p_i(n)\right\rceil: p_i(n) \in \widehat{SGP_d} \text{ and } c_i \in \mathbb{Z}\}.$$
	The  {\it special integer-valued generalized polynomials}  is then defined by
	$$\widetilde{SGP}=\bigcup_{d=1}^{\infty} \widetilde{SGP_d}.$$
\end{defn}
Clearly $\widetilde{SGP} \subset \mathcal{G}$ and we have the following obsevation.

\begin{lem}\label{lem-sum-sgp}
	Let $p_1, \cdots, p_d\in \widehat{SGP_s}$ (for some $s\in \mathbb{N}$).
	Then for any $n\in \mathbb{Z}$ with
	$$-\frac{1}{2}<\{p_1(n)\}+\cdots + \{p_d(n)\}< \frac{1}{2},$$
	we have $\left\lceil p_1(n)+\cdots + p_d(n) \right\rceil= \sum\limits_{i=1}^{d}  \left\lceil p_i(n)\right\rceil$.
\end{lem}	

The following lemma shows the the relationship 
between $\widetilde{GP_d}$ and $\widetilde{SGP_d}$.
\begin{lem}\label{lem sgp_gp}
	Let $d \in \mathbb{N}$ and $p(n) \in \widetilde{GP_d}$. Then there exist $h(n)\in \widetilde{SGP_d}$ and a set
	$$C=C(\delta,q_1,\cdots,q_t)=\bigcap_{k=1}^{t} \{n\in  \mathbb{Z}: \{q_k(n)\}\in (-\delta,\delta)\}$$
	such that $$p(n)=h(n), \forall n \in C,$$ where $\delta>0$ is small enough and
	$q_k \in \widehat{SGP_d}, k=1,2,\dots, t$ for some $t \in \mathbb{N}$.
	
\end{lem}
\begin{proof}
	We will prove it by induction on $d$.
	
	When $d =1$, we may assume that 
	$p(n)= \left\lceil \sum_{j=1}^{m} \alpha_j \left \lceil \beta_j n \right \rceil \right \rceil $. Let
	$$q_j(n)=\alpha_j \left \lceil \beta_j n \right \rceil, j=1,\dots,m.$$
	Let $0<\delta<\frac{1}{2m}$, we set
	$$C=C(\delta,q_1,\dots,q_m)=\bigcap_{j=1}^{m} \{n \in \mathbb{Z}: \{q_j(n)\} \in (-\delta,\delta)\}.$$ 
	Since for each $n \in C$, $\{q_j(n)\}=\{\alpha_j \left \lceil \beta_j n \right \rceil \} \in (-\delta,\delta)$,  
	$$-\frac{1}{2}<-m\delta<\sum_{j=1}^{m}\{\alpha_j \left \lceil \beta_j n \right \rceil \}<m\delta<\frac{1}{2}, \forall n \in C.$$
	Let $h(n)=\sum_{j=1}^{m}\left\lceil  \alpha_j \left \lceil \beta_j n \right \rceil \right \rceil$, then $h(n) \in \widetilde{SGP_1}$.
	Hence by Lemma \ref{lem-sum-sgp}, $p(n)= h(n), \forall n \in C.$
	
	Assmume that the result holds for $d>1$. Next we will show the result holds for $d+1$. 
	We just need to show that when $p(n)= \left\lceil r(n)\right\rceil$ the result holds, where
	$$r(n)=a_0n^{p_0}\left\lceil {f_1(n) } \right\rceil \cdots \left\lceil {f_k(n)} \right\rceil $$
	(with $a_0\in \mathbb{R}, p_0 \ge 0, k\ge 0, f_l\in GP_{p_l}$ and $\sum_{l=0}^{k}p_l=d+1$).
	
	If $p_0={d+1}$, then $p(n)=\left\lceil a_0n^{d+1}\right\rceil \in \widetilde{SGP_{d+1}}$.
	Next we assume that $0\le p_0< d+1$ and $0<p_l <d+1, l=1,2,\dots,k$. For each $1 \le l \le k$, by induction hypothesis, there exist $h_l(n) \in \widetilde{SGP_{p_l}} $ and 
	$C_l$ such that 
	$$\left\lceil {f_l(n) } \right\rceil =h_l(n):=\sum_{i=1}^{b_l}c_{l,i}\left\lceil r_{l,i}(n) \right\rceil, \forall n \in C_l,$$
	where $c_{l,i} \in \mathbb{Z}$, $r_{l,i}(n)\in \widehat{SGP_{p_l}}$    and 
	$$C_l=C_l(\delta_l, q_{l,1},\dots, q_{l,t_l})$$
	(with $\delta_l>0$ is small enough and $q_{l,k} \in \widehat{SGP_{p_l}}, k=1, \cdots, t_l$ for some $t_l \in \mathbb{N}$). 
	
	For any $n \in \bigcap_{l=1}^{k}C_l$, 
	\begin{eqnarray*}
		r(n)&=&a_0n^{p_0}\left\lceil {f_1(n) } \right\rceil \cdots \left\lceil {f_k(n)} \right\rceil\\
		&=&a_0n^{p_0}h_1(n)\cdots h_k(n)\\
		&=&a_0n^{p_0}(\sum_{i=1}^{b_1}c_{1,i}\left\lceil r_{1,i}(n) \right\rceil)\cdots
		(\sum_{i=1}^{b_k}c_{k,i}\left\lceil r_{k,i}(n) \right\rceil).
	\end{eqnarray*}
	Note that $\left\lceil r_{1,i_1}(n) \right\rceil \cdots \left\lceil r_{k,i_k}(n) \right\rceil \in \widehat{SGP_{d+1-p_0}} $ and are integer-valued, then $r(n)$ can be written as
	$$r(n)=\sum_{j=1}^{m}\beta_jn^{p_0}\left\lceil d_j(n)\right\rceil, $$
	where $d_j(n)$ is of the form $\left\lceil r_{1,j_1}(n) \right\rceil \cdots \left\lceil r_{k,j_k}(n) \right\rceil $.
	
	Let $Q =\{ \beta_1 n^{p_0}\left\lceil d_1(n)\right\rceil,\dots, \beta_m n^{p_0}\left\lceil d_m(n)\right\rceil \} \cup (\bigcup_{l=1}^{k}\{q_{l,1}(n),\dots,q_{l,t_l}(n)\})$. Let
	$0<\delta<\min \{\frac{1}{2m}, \delta_1,\dots,\delta_k\}$
	and $$C=C(\delta,Q)=\bigcap_{q(n) \in Q}\{n \in \mathbb{Z}:\{q(n)\}\in (-\delta,\delta) \}.$$
	Clearly $C\subset \bigcap_{l=1}^{k}C_l$. For each $n \in C$,
	$\{\beta_jn^{p_0} \left\lceil d_j(n)\right\rceil\}\in (-\delta,\delta),j=1,2,\dots,m.$ Hence
	$$-\frac{1}{2}<-m\delta<\sum_{j=1}^m\{\beta_jn^{p_0} \left\lceil d_j(n)\right\rceil\}<m\delta<\frac{1}{2}.$$
	Let $h(n)=\sum_{j=1}^m\left\lceil \beta_jn^{p_0} \left\lceil d_j(n)\right\rceil\right\rceil$, $h(n)\in \widetilde{SGP_{d+1}}$. By Lemma \ref{lem-sum-sgp}, $p(n)=h(n), \forall n\in C$.

\end{proof}

\medskip

By Lemma \ref{lem sgp_gp}, every $p(n) \in \widetilde{GP_d}$ correspondes to an $h(n) \in \widetilde{SGP_d}$,  we call the maximal-degree components of $h(n)$ be the maximal-degree components of $p(n)$. But we need to mention that here we will not do the $+$ and $-$. 
For instance, let 
$p(n)=n\left\lceil 2 \pi n^2 - \left\lceil 2 \pi n^2 \right\rceil +\sqrt{2}n\right \rceil $ then we
denote $h(n):=n\left \lceil 2\pi n^2 \right\rceil-n\left \lceil 2\pi n^2 \right\rceil +n\left\lceil \sqrt{2}n \right\rceil$,
and we denote the maximal-degree components of $p(n),h(n)$ be $n\left \lceil 2\pi n^2 \right\rceil$ and $-n\left \lceil 2\pi n^2 \right\rceil $ and the coefficients of the maximal-degree components of $p(n), h(n)$ are $2\pi$ and $-2\pi$.

\begin{defn}\label{def-nondegenerate}	
	Let $p(n) \in \mathcal{G}$, we denote $A(p(n))$ be the sum of the coefficients of the maximal-degree componentes of $p(n)$.  
	Let $p_1,p_2,\cdots,p_d \in  \mathcal{G}$, a tuple $(p_1,p_2,\cdots,p_d)$ is called a
	{\it non-degenerate tuple} if $A(p_i) \neq 0$ and $A(p_i-p_j) \neq 0$, $1 \le i\neq j \le d$. 
\end{defn}

For instance, 
$A(\left\lceil an^2 \left \lceil bn\right\rceil +\left\lceil cn^3\right\rceil \right\rceil +dn^3 +2n^2)=ab+c+d$, $A(n+n\left \lceil 2\pi n -\left\lceil 2\pi n\right\rceil \right\rceil)=0$.	
$(n^2+n,n^2+\left\lceil \sqrt{3}n\right\rceil)$ is non-degenerate, $(n\left\lceil 2\pi n \right\rceil+n, \left\lceil 2\pi n^2 \right\rceil+2n)$ is not non-degenerate.

\medskip

The key ingredient in the proof of the main result is to view the integer-valued generalized polynomials, in some sense, as the ordinary polynomials. To do this, we need to introduce the following definition.
\begin{defn}\label{def proper}
	Let $p(n)\in \widetilde{SGP}$, $m \in \mathbb{Z}$ and $C \subset \mathbb{Z}$. We say that $p$ is \emph{proper} with respect to (w.r.t. for short) $m$ and $C$ if
	for every $n\in C$,
	\begin{itemize}
		\item if $deg(p)=1$, $p(n+m)=p(n)+p(m)$.
		\item if $deg(p)>1$, 	$p(n+m)-p(n)-p(m)=q(n),$
		where $q(n)\in \widetilde{SGP}$ and $deg(q) <deg(p)$.
	\end{itemize}
\end{defn}

For example, let $p(n)=\left\lceil an^2\right\rceil$, if
$$p(n+m)=\left\lceil a(n+m)^2\right\rceil= \left\lceil an^2\right\rceil+\left\lceil am^2\right\rceil+\left\lceil 2amn\right\rceil, \forall n \in C,$$ then we say $p(n)$ is proper w.r.t. $m$ and $C$.

Let $p(n) \in \widetilde{SGP}$, $m\in \mathbb{Z}$. To study whether there exists $C$ such that $p(n)$ is proper w.r.t. $m$ and $C$,  we need to introduce the following notion.
\begin{defn}\label{def-m-good}
	Let $p(n) \in \widetilde{SGP}$ and $m\in \mathbb{Z}$. 
	\begin{itemize}
		\item If $p(n)=\left\lceil L(a_1n^{j_1},\dots, a_ln^{j_l}) \right \rceil$, we say $m$ {\it is good w.r.t.} $p(n)$ if for any $1\le t \le l$,  $\{L(a_{t}m^{j_{t}},a_{t+1}m^{j_{t+1}},\cdots,a_lm^{j_l})\} \neq \frac{1}{2} $.
		\item If $p(n)=\left\lceil \prod_{i=1}^{k}r_{i}(n) \right\rceil$
		with $r_i(n)=L(a_{i,1}n^{j_{i,1}},\dots, a_{i,l_i}n^{j_{i,l_i}})$, we say $m$ is good w.r.t. $p(n)$ if	
		$\{\prod_{i=1}^{k}r_i(m)\} \neq \frac{1}{2}$ 	
		and 	 $m$ is
		good w.r.t.  $\left\lceil r_i(n) \right \rceil$ for each $1 \le i \le k$.
		\item	If $p(n)=\sum_{t=1}^{k}c_t\left\lceil q_t(n) \right\rceil$ with $c_i \in\mathbb{Z}$ and each $q_t(n)$ is of the form $\prod_{i=1}^{k}r_{i}(n)$ with $r_i(n)=L(a_{i,1}n^{j_{i,1}},\dots, a_{i,l_i}n^{j_{i,l_i}})$, we say $m$ is good
		w.r.t. $p(n)$ if $m$ is good w.r.t. $\left\lceil q_t(n)\right\rceil$ for each $1 \le t \le k$. 
		
	\end{itemize}
\end{defn}
For example, if $\{ bm\left\lceil cm \right\rceil\} \neq \frac{1}{2}$ and $\{ cm \} \neq \frac{1}{2}$,   then $m$ is good w.r.t. $p(n)=\left\lceil bn\left\lceil cn \right\rceil\right\rceil$.

We have the following observation.
\begin{lem}\label{lem-goodm}
	Let $p(n) \in \widetilde{SGP}$. Then there exist $\delta>0$, $Q\subset \widehat{SGP_s}$ (for some $s\in \mathbb{N}$) and
	$$C(\delta,Q)=\bigcap_{q(n)\in Q}\{n\in \mathbb{Z}: \{q(n)\} \in (-\delta,\delta)\}$$
	such that for each $m\in C(\delta,Q)$, $m$ is good w.r.t. $p(n)$.
\end{lem}
\begin{proof}
	Choose $0<\delta<\frac{1}{4}$. 
	\begin{itemize}
		\item If $p(n)=\left\lceil L(a_1n^{j_1},\dots, a_ln^{j_l}) \right \rceil$, let $Q=\{L(a_{t}n^{j_{t}},a_{t+1}n^{j_{t+1}},\cdots,a_ln^{j_l}): 1\le t \le l \}$. Then for each $m\in C(\delta,Q)$, $m$ is good w.r.t. $p(n)$.
		\item If $p(n)=\left\lceil \prod_{i=1}^{k}r_{i}(n) \right\rceil$
		with $r_i(n)=L(a_{i,1}n^{j_{i,1}},\dots, a_{i,l_i}n^{j_{i,l_i}})$, let
		$$Q=\{\prod_{i=1}^{k}r_{i}(n)\}
		\cup\bigcup_{i=1}^k\{L(a_{i,t}n^{j_{i,t}},a_{i,t+1}n^{j_{t+1}},\cdots,a_{i,l_i}n^{j_{i,l_i}}): 1\le t \le l \}.$$
		Then for each $m \in C(\delta, Q)$, $m$ is good w.r.t. $p(n)$.
		\item If $p(n)=\sum_{t=1}^{k}c_t\left\lceil q_t(n) \right\rceil$. For each
		$\left\lceil q_t(n) \right\rceil$, by the above argument there exists a $Q_t$ such that for each $m\in C(\delta,Q_t)$, $m$ is good w.r.t. $\left \lceil q_t(n) \right\rceil$.
		Let $Q=(\bigcup_{t=1}^{k}Q_t)$, for each $m\in C(\delta,Q)$, $m$ is good w.r.t. $p(n)$.
	\end{itemize}
	
\end{proof}

\medskip
The following  lemmas  are very useful in our proof. 
We first prove the simple case  to illustrate our idea.
The general case  can be deduced directly.

\begin{lem}\label{lem-proper-simple}
	Let $p(n)=\left\lceil r(n) \right\rceil, n\in \mathbb{Z}$, where $r(n)\in \widehat{SGP_d}$ for some $d \in \mathbb{N}$. 
	Let $l\in \mathbb{N}$, $m_i \in \mathbb{Z}$ and $m_i$ is good w.r.t. $p(n)$ for each $1 \le i \le l$.
	Then for any $\varepsilon >0$,
	there exists
	$$C=C(\delta, q_1, \cdots, q_t)=\bigcap_{k=1}^{t}\{n\in \mathbb{Z}: \{q_k(n)\}\in (-\delta,\delta)\},$$
	where  $\delta>0 \ (\delta<\varepsilon)$ is a small enough number, 	
	and $q_k \in \widehat{SGP_d}, k=1,2,\dots,t$ for some $t\in \mathbb{N}$,
	such that for all $ i\in\{1, \cdots, l\},$
	\begin{enumerate}
		\item $p(n)$ is proper w.r.t. $m_i$ and $C$.
		\item $\{r(n+m_i)\} \in (\{r(m_i)\}-\varepsilon,\{r(m_i)\}+\varepsilon),\forall n \in C$.
	\end{enumerate}
\end{lem}

\begin{proof}
	We first show a special case $r(n)=bn\left\lceil cn \right\rceil$ to illustrate our idea, the general cases are similar.
	
	Let $\delta_1=\frac{1}{2}- \mathop {\max }_{i=1,\dots,l} \{|\{bm_i\left\lceil cm_i\right\rceil \}|, |\{cm_i\}|\}  $. Since
	for each $1\le i\le l$, $m_i$ is good w.r.t. $\left\lceil r(n)\right\rceil$, $\delta_1>0$.
	Choose $0<\delta<\min \{\frac{\delta_1}{4},\frac{\varepsilon}{3}\}$ and let
	\begin{eqnarray*}
		C(\delta)=\bigcap_{i=1}^{l}\{n \in \mathbb{Z}: \{bn\left\lceil cn\right\rceil \}, \{bn\left\lceil cm_i\right\rceil \}, \{bm_i\left\lceil cn\right\rceil \}, \{cn\} \in (-\delta,\delta)\}.
	\end{eqnarray*}
	
	Since for all $ i=1, \cdots, l$ and $n\in C(\delta)$, we have
	$$|\{cm_i\}| \le \frac{1}{2}-\delta_1,  \{cn\}\in (-\delta, \delta),$$
	$$|\{bm_i\left\lceil cm_i\right\rceil \}|\le \frac{1}{2}-\delta_1,
	\{bn \left\lceil cm_i\right\rceil \}, \{bn \left\lceil cn\right\rceil \}, \{bm_i \left\lceil cn\right\rceil \}\in (-\delta,\delta).$$
	Then 
	\begin{equation}\label{eq-lem-prop-1}
	-\frac{1}{2}< \{cm_i\}+\{cn\}<\frac{1}{2},
	\end{equation}
	\begin{equation}\label{eq-lem-prop-2}
	-\frac{1}{2}<\{bn \left\lceil cn\right\rceil \}+\{bn \left\lceil c m_i\right\rceil \} +\{bm_i \left\lceil cn\right\rceil \}+ \{bm_i\left\lceil cm_i\right\rceil \}<\frac{1}{2}.
	\end{equation}
	By \eqref{eq-lem-prop-1} and Lemma \ref{lem-sum-sgp}, $\left\lceil cm_i+cn\right\rceil=\left\lceil cm_i\right\rceil+\left\lceil cn\right\rceil$.
	Then 
	\begin{eqnarray*}
		r(n+m_i)&=&b(n+m_i)\left\lceil c(n+m_i)\right\rceil \\
		&=&(bn+bm_i)(\left\lceil cn\right\rceil+\left\lceil cm_i\right\rceil)\\
		&=& bn\left\lceil cn\right\rceil+bn\left\lceil cm_i \right\rceil+bm_i\left\lceil cm_i\right\rceil+bm_i\left\lceil cm_i\right\rceil.
	\end{eqnarray*}
	By \eqref{eq-lem-prop-2} and Lemma \ref{lem-sum-sgp},
	\begin{eqnarray*}
		p(n+m_i)&=&\left\lceil r(n+m_i)\right\rceil \\
		&=& \left\lceil bn\left\lceil cn\right\rceil\right\rceil+\left\lceil bn\left\lceil cm_i \right\rceil\right\rceil+ \left\lceil bm_i\left\lceil cn\right\rceil\right\rceil+\left\lceil bm_i\left\lceil cm_i\right\rceil\right\rceil\\ 
		&=& p(n)+p(m_i)+(\left\lceil bn\left\lceil cm_i \right\rceil\right\rceil+\left\lceil bm_i\left\lceil cn\right\rceil\right\rceil).
	\end{eqnarray*}
	Which implies $p(n+m_i)$ is proper.
	It also implies that
	\begin{eqnarray*}
		\{r(n+m_i)\}&=&\{r(m_i) +bn \left\lceil c n\right\rceil +bn \left\lceil c m_i\right\rceil +b m_i \left\lceil c n \right\rceil\}\\
		&\in& (\{r(m_i) \}-\varepsilon,\{r(m_i) \}+\varepsilon).
	\end{eqnarray*}

	We will prove the general cases by proving the result holds for the following three cases.
	
	\noindent {\bf Case 1: $r(n)$ is of the form $r(n)=L(an^j)$.} 
	
	For any $1 \le i \le l$. Let
	$\tilde{Q}_i$ be the set of all the components of the expansion of $a(n+m_i)^j$  and $$Q_i=\tilde{Q}_i\setminus \{L(am_i^j)\}$$
	(e.g.  $a(n+m_i)^2=an^2+2anm_i+am_i^2$, $\tilde{Q}_{i}=\{an^2, 2anm_i, am_i^2\}$ and  $Q_i=\{an^2,2anm_i\}$). It is clear that $\#\tilde{Q}_i \le 2^j$, where $\# Q$ is the number of elements of the set $Q$.
	
	Let $\delta_1=\frac{1}{2}-\max_{i=1,2,\cdots,l}\{|{L(am_i^j)}|\}$. Since for each $i=1, \cdots, l$, $m_i$ is good w.r.t. $p(n)=\left\lceil r(n) \right\rceil$, $\delta_1 >0$. Choose 
	$0<\delta<\min \{\frac{\delta_1}{2^j},\frac{\varepsilon}{2^j-1}\}$ and let
	$$C(\delta)=\bigcap_{i=1}^{l}\bigcap_{q(n) \in Q_i }\{n\in \mathbb{Z}: \{q(n) \} \in (-\delta,\delta)\}.$$
	For any $1\le i \le l$ and any $n \in C(\delta)$, 
	$$-\delta_1<-(2^j-1)\delta<\sum_{q(n) \in Q_i}\{q(n)\} <(2^j-1)\delta<\delta_1,$$
	$$|\{L(am_i^j)\}| \le \max_{i=1,2,\cdots,l}\{|{L(am_i^j)}|\}=\frac{1}{2}-\delta_1.$$
	Since
	\begin{eqnarray*}
		r(n+m_i)=a(n+m_i)^j=\sum_{q(n) \in Q_i} q(n)+L(am_i^j)=\sum_{q(n) \in Q_i} q(n)+r(m_i)
	\end{eqnarray*}
	and 
	\begin{eqnarray*}
		-\frac{1}{2}=-\delta_1-(\frac{1}{2}-\delta_1)	<\sum_{q(n) \in Q_i} \{q(n)\}+\{L(am_i^j)\}<\delta_1+\frac{1}{2}-\delta_1=\frac{1}{2}.
	\end{eqnarray*}
	By Lemma \ref{lem-sum-sgp},
	\begin{eqnarray*}
		p(n+m_i)&=&\left\lceil r(n+m_i) \right\rceil =\left\lceil \sum_{q(n) \in Q_i} q(n)+L(am_i^j) \right\rceil=\sum_{q(n) \in Q_i} \left\lceil  q(n) \right\rceil+\left\lceil L(am_i^j) \right\rceil \\
		&=& p(n)+p(m_i)+\sum_{q(n) \in Q_i\setminus\{r(n)\}} \left\lceil  q(n) \right\rceil.
	\end{eqnarray*}
	Which implies $p(n)$ is proper w.r.t. $m_i$ and $C$.
	Since 
	$$-\varepsilon<-(2^j-1)\delta<\sum_{q(n) \in Q_i}\{q(n)\} <(2^j-1)\delta<\varepsilon.$$
	We have
	$$\{r(n+m_i)\}=\{r(m_i)+\sum_{q(n) \in Q_i} q(n)\}=\{r(m_i)\}+\{\sum_{q(n) \in Q_i} q(n)\} \in (\{r(m_i)\}-\varepsilon,\{r(m_i)\}+\varepsilon).$$
	
	\noindent {\bf Case 2: $r(n)$ is of the form $r(n)=L(a_1n^{j_1},\dots,a_tn^{j_t})$.}
	
	For any $1 \le  i \le l$ and $1 \le k \le t$,
	let $\tilde{Q}_{i,k}$ be the set of all the components of $a_k(n+m_i)^{j_k}$, we denote 
	\begin{eqnarray*}
		&&Q_{i,t}=\tilde{Q}_{i,t} \setminus \{L(a_tm_i^{j_t})\}, \\ &&Q_{i,t-1}=\tilde{Q}_{i,t-1}\left\lceil\tilde{Q}_{i,t}\right\rceil \setminus \{L(a_{t-1}m_i^{j_{t-1}},a_tm_i^{j_t})\},\\
		&&\cdots\cdots,\\
		&&Q_{i,1}=\tilde{Q}_{i,1}\left\lceil \tilde{Q}_{i,2} \left\lceil \cdots \left\lceil \tilde{Q}_{i,t} \right\rceil \cdots \right\rceil \right\rceil \setminus \{L(a_{1}m_i^{j_{1}},\cdots,a_tm_i^{j_t})\},
	\end{eqnarray*}
	and $$Q_i=Q_{i,t}\cup Q_{i,t-1}\cup \cdots \cup Q_{i,1},$$
	where $\left\lceil A\right\rceil:=\{\left\lceil a \right\rceil: a \in A\}$ and $AB:=\{ab:a \in A,b\in B\}$ for $A,B \subset \mathcal{G}$.

	Let $$\delta_1=\frac{1}{2}-\max_{i=1,2,\dots,l}\{|\{L(a_tm_i^{j_t})\}|,
	|\{L(a_{t-1}m_i^{j_{t-1}},a_tm_i^{j_t})\}|,\cdots,|\{L(a_{1}m_i^{j_{1}},\cdots,a_tm_i^{j_t})\}|\}.$$
	Since $m_i$ is good w.r.t. $p(n)=\left\lceil r(n) \right\rceil$, $\delta_1 >0$.
	Let $$L:=2^{j_t}+2^{j_t+j_{t-1}}+\dots+2^{j_t+j_{t-1}+\dots+j_1}>\# Q_i,$$
	we  choose 
	$0<\delta<\min \{\frac{\delta_1}{L},\frac{\varepsilon}{L-1}\}$. Let 
	$$C(\delta)=\bigcap_{i=1}^{l}\bigcap_{q(n) \in Q_i }\{n\in \mathbb{Z}: \{q(n) \} \in (-\delta,\delta)\}.$$
	For any $1 \le i \le l$ and any $n \in C(\delta)$,
	$$-\delta_1<-2^{j_t}\delta<\sum_{q(n) \in Q_{i,t}}\{q(n)\}<2^{j_t}\delta<\delta_1,$$
	$$|\{L(a_tm_i^{j_t})\}| \le \frac{1}{2}-\delta_1,$$
	using the same argument as in case 1 and applying Lemma \ref{lem-sum-sgp}, we have
	$$\left\lceil L(a_t(n+m_i)^{j_t})\right\rceil=\left\lceil L(a_tm_i^{j_t})\right\rceil+\sum_{q\in Q_{i,t}}\left\lceil q(n) \right\rceil.$$
	Then 
	\begin{eqnarray*}
		L(a_{t-1}(n+m_i)^{j_{t-1}},a_t(n+m_i)^{j_t})&=&( \sum_{q\in \tilde{Q}_{i,t-1}} q(n) )(\left\lceil L(a_tm_i^{j_t})\right\rceil+\sum_{q\in Q_{i,t}}\left\lceil q(n) \right\rceil)\\
		&=&  L(a_{t-1}m_i^{j_{t-1}},a_tm_i^{j_t})+\sum_{q\in Q_{i,t-1}} q(n).
	\end{eqnarray*}
	Since
	$$-\delta_1<-2^{j_t+j_{t-1}}\delta<\sum_{q(n) \in Q_{i,t-1}}\{q(n)\}<2^{j_t+j_{t-1}}\delta<\delta_1,$$
	$$|\{L(a_{t-1}m_i^{j_{t-1}},a_tm_i^{j_t})\}| \le \frac{1}{2}-\delta_1,$$
	using the same argument as in case 1 and applying Lemma \ref{lem-sum-sgp}, we have
	\begin{eqnarray*}
		\left\lceil L(a_{t-1}(n+m_i)^{j_{t-1}},a_t(n+m_i)^{j_t})\right\rceil
		= \left\lceil L(a_{t-1}m_i^{j_{t-1}},a_tm_i^{j_t}) \right\rceil +\sum_{q\in Q_{i,t-1}} \left\lceil q(n) \right\rceil.
	\end{eqnarray*}
	Inductively, we have
	\begin{eqnarray*}
		L(a_{1}(n+m_i)^{j_{1}},\cdots,a_t(n+m_i)^{j_t})&=&( \sum_{q\in \tilde{Q}_{i,1}} q(n) )(\left\lceil L(a_{2}m_i^{j_{2}},\cdots,a_tm_i^{j_t})\right\rceil+\sum_{q\in Q_{i,2}}\left\lceil q(n) \right\rceil)\\
		&=&L(a_{1}m_i^{j_{1}},\cdots,a_tm_i^{j_t})+\sum_{q(n)\in Q_{i,1}}q(n)  .
	\end{eqnarray*}
	Since
	$$-\delta_1<-2^{j_t+j_{t-1}+\cdots+j_1}\delta<\sum_{q(n) \in Q_{i,1}}\{q(n)\}<2^{j_t+j_{t-1}+\cdots+j_1}\delta<\delta_1,$$
	$$|\{L(a_{1}m_i^{j_{1}},\cdots,a_tm_i^{j_t})\}| \le \frac{1}{2}-\delta_,$$
	using the same argument as in case 1, we have
	$$-\frac{1}{2}<\{L(a_{1}m_i^{j_{1}},\cdots,a_tm_i^{j_t})\}+\sum_{q\in Q_{i,1}} \{q(n)\}<\frac{1}{2}.$$
	Then applying Lemma \ref{lem-sum-sgp}, 
	\begin{eqnarray*}
		\left\lceil r(n+m_i)  \right\rceil &=&\left\lceil L(a_{1}(n+m_i)^{j_{1}},\cdots,a_t(n+m_i)^{j_t})\right\rceil
		= \left\lceil L(a_{1}m_i^{j_{1}},\cdots,a_tm_i^{j_t}) \right\rceil +\sum_{q\in Q_{i,1}} \left\lceil q(n) \right\rceil\\
		&=&\left\lceil r(m_i)  \right\rceil+\left\lceil r(n)  \right\rceil+\sum_{q\in Q_{i,1}\setminus \{r(n)\}} \left\lceil q(n) \right\rceil.
	\end{eqnarray*}
	It implies $p(n)=\left\lceil r(n) \right\rceil$ is proper w.r.t. $m_i$ and $C(\delta)$.
	Since
	$$-\varepsilon<-(L-1)\delta<\sum_{q(n) \in Q_{i,1}}\{q(n)\}<(L-1)\delta<\varepsilon,$$
	we have 
	\begin{eqnarray*}
		\{r(n+m_i)\}&=&\{ L(a_{1}m_i^{j_{1}},\cdots,a_tm_i^{j_t})
		+\sum_{q\in Q_{i,1}}  q(n)  \}\\
		&\in& (\{r(m_i)\}-\varepsilon,\{r(m_i)\}-\varepsilon).
	\end{eqnarray*}
	
	\noindent{\bf Case 3: $r(n)$ is of the form $r(n)=\prod_{h=1}^{k}r_h(n)$ with $r_h(n)=L(a_{h,1}n^{j_{h,1}},\dots,a_{h,t_h}n^{j_{h,t_h}})$.}
	
	Notice that 
	$$L(an,bn)L(cn,dn)=(an\left\lceil bn \right \rceil)(cn\left\lceil dn \right \rceil)=acn^2\left\lceil bn \right \rceil\left\lceil dn \right \rceil=L(acn^2,bn)\left\lceil L(dn)\right\rceil,$$ we can assume 
	$r(n)=r_1(n)\prod_{h=2}^{k}\left\lceil r_h(n)\right\rceil$ with $r_h(n)=L(a_{h,1}n^{j_{h,1}},\dots,a_{h,t_h}n^{j_{h,t_h}}), h\ge 1$.
	
	For each $1 \le h \le k $ and $1 \le i \le l$, by case 2, there exsit $L_h \in \mathbb{Z}$, $\delta_h>0$, $  \tilde{Q}^h_{i,1},\cdots, \tilde{Q}^h_{i,t_h}, Q^h_{i,1},\cdots, Q^h_{i,t_h}$ with
	$$C(\delta_h)=\bigcap_{i=1}^{l}\bigcap_{q(n)\in Q_i^h}\{n\in \mathbb{Z}:\{q(n)\}\in (-\delta_h,\delta_h)\},$$
	$$Q_i^h=Q^h_{i,t_h}\cup Q^h_{i,t_h-1}\cup \cdots \cup Q^h_{i,1},$$ 
	corresponding to $r_h(n)$ and $m_i$, such that
	$$\left\lceil r_h(n+m_i)\right\rceil= \left\lceil r_h(m_i)\right\rceil+\sum_{q(n)\in Q^h_{i,1}} \left\lceil q(n)\right\rceil, \forall n \in C(\delta_h). $$
	
	Let $Q_i=\tilde{Q}_{i,1}^1\prod_{h=2}^{k}\left\lceil \tilde{Q}_{i,1}^h\right\rceil \backslash \{r(m_i)\}$, $L=\prod_{h=1}^{k}L_h$ and $\tilde{\delta}=\min_{1 \le h \le k}\delta_h$.
	Let 
	\begin{eqnarray*}
		B=&& \bigcup_{i=1}^{l}\bigcup_{h=1}^{k}\{L(a_{h,t_h}m_i^{j_{h,t_h}}),
		L(a_{h,t_h-1}m_i^{j_{h,t_h-1}},a_{h,t_h}m_i^{j_{h,t_h}}),\cdots,L(a_{h,1}m_i^{j_{h,1}},\cdots,a_{h,t_h}m_i^{j_{h,t_h}})\} \\
		&&\bigcup \{r_1(m_i)\prod_{h=2}^{k}\left\lceil r_h(m_i)\right\rceil: i=1,2,\dots,l\}
	\end{eqnarray*}
	and $\hat{\delta}=\frac{1}{2}-\max_{q\in B}\{|\{q\}|\}$. Since $m_i$ is good 
	w.r.t. $p(n)$, $\hat{\delta}>0$. We choose $\delta< \min\{\frac{\varepsilon}{L},\frac{\tilde{\delta}}{L},\frac{\hat{\delta}}{L}\}$ and let
	$$C(\delta)=(\bigcap_{h=1}^{k}\bigcap_{i=1}^{l}\bigcap_{q(n)\in Q_i^h}\{n\in \mathbb{Z}:\{q(n)\}\in (-\delta,\delta)\}) \cap (\bigcap_{i=1}^{l}\bigcap_{q(n)\in Q_i}\{n\in \mathbb{Z}:\{q(n)\}\in (-\delta,\delta)\}).$$
	Since $C(\delta) \subset \bigcap_{h=1}^{k}C(\delta_h)$, for any $m_i$ and $n\in C(\delta)$,
	\begin{eqnarray*}
		r(n+m_i)
		&=&r_1(n+m_i)\prod_{h=2}^{k}\left\lceil r_h(n+m_i) \right\rceil\\
		&=&(r_1(m_i)+\sum_{q(n) \in Q_{i,1}^{1}}\left \lceil q(n)\right\rceil)
		\prod_{h=2}^{k}(\left \lceil r_h(m_i)\right\rceil+\sum_{q(n) \in Q_{i,1}^{h}}\left \lceil q(n)\right\rceil)\\
		&=& r_1(m_i)\prod_{h=2}^{k}\left \lceil r_h(m_i)\right\rceil+\sum_{q(n) \in Q_{i}} q(n)
	\end{eqnarray*}
	and
	$$-\hat{\delta}<-L\delta<\sum_{q(n)\in Q_i}\{q(n)\}<L\delta<\hat{\delta}, \ |\{ r_1(m_i)\prod_{h=2}^{k}\left \lceil r_h(m_i)\right\rceil  \}|< \frac{1}{2}-\hat{\delta}.$$
	Hence 
	$$-\frac{1}{2}<\{r_1(m_i)\prod_{h=2}^{k}\left \lceil r_h(m_i)\right\rceil\}+\sum_{q(n)\in Q_i}\{q(n)\}<\frac{1}{2}.$$
	By Lemma \ref{lem-sum-sgp}, 
	\begin{eqnarray*}
		\left\lceil r(n+m_i) \right\rceil
		&=& \left \lceil r_1(m_i)\prod_{h=2}^{k}\left \lceil r_h(m_i)\right\rceil \right\rceil+\sum_{q(n) \in Q_{i}} \left \lceil  q(n)\right\rceil\\
		&=&  \left \lceil r(m_i)\right\rceil + \left \lceil r(n)\right\rceil +  \sum_{q(n) \in Q_{i}\backslash \{r(n)\}} \left \lceil  q(n)\right\rceil.
	\end{eqnarray*}
	It implies $p(n)$ is proper w.r.t. $m_i$ and $C(\delta)$.
	Since 
	$$ -\varepsilon<-L\delta<\sum_{q(n)\in Q_i}\{q(n)\}<L\delta<\varepsilon,$$
	we have
	\begin{eqnarray*}
		\{r(n+m_i)\}
		&=& \{ r_1(m_i)\prod_{h=2}^{k}\left \lceil r_h(m_i)\right\rceil+\sum_{q(n) \in Q_{i}} q(n)\} \\
		&\in& (\{ r_1(m_i)\prod_{h=2}^{k}\left \lceil r_h(m_i)\right\rceil\}-\varepsilon, \{r_1(m_i)\prod_{h=2}^{k}\left \lceil r_h(m_i)\right\rceil\}+\varepsilon)\\
		&=& ( \{r(m_i)\} -\varepsilon, \{r(m_i)\}+\varepsilon).
	\end{eqnarray*}
	
	Thus we finish the proof.
\end{proof}


\medskip
Since $\mathcal{F}_{d, 0}$ is a filter, we have the following result.

\begin{lem} \label{asso for special cases}
	Let  $p_1(n)=\left\lceil r_1(n) \right\rceil, \cdots, p_t(n)=\left\lceil r_t(n) \right\rceil, n\in \mathbb{Z}$, where $r_i \in \widehat{SGP_d}, i=1, \cdots, t$ for some $d \in \mathbb{N}$. Let $l\in \mathbb{N}, m_j \in \mathbb{Z}$ and  $m_j$ is good w.r.t.
	$p_i(n)$ for $ 1\le j \le l, 1\le i \le t$.
	Then for any $\varepsilon>0$,
	there exists
	$$C=C(\delta)=\bigcap_{k=1}^{h}\{n\in \mathbb{Z}: \{q_k(n)\}\in (-\delta,\delta)\},$$
	where  $\delta>0 \ (\delta< \varepsilon)$ is a small enough number, 	
	$s = \max_{1\le i\le t}deg (p_i)$ and $q_k \in \widehat{SGP_s}, k=1,2,\dots,h$ for some $h\in \mathbb{N}$,
	such that for all $i\in \{1, \cdots, t\}, j\in\{1, \cdots, l\}$,
	\begin{enumerate}
		\item $p_i(n)$ is proper w.r.t. $m_j$ and $C$.
		\item $\{r_i(n+m_j)\} \in (\{r_i(m_j)\}-\varepsilon,\{r_i(m_j)\}+\varepsilon), \forall n \in C$.
	\end{enumerate}
\end{lem}

And the  general case is the following lemma.

\begin{lem} \label{lem-proper}
	Let $p_1, \cdots, p_d\in \widetilde{SGP}$. Let $l\in \mathbb{N}, m_j\in \mathbb{Z}$ and $m_j$ is good w.r.t. $p_i(n)$ for $1\le i \le d, 1 \le j \le l$.  
	Then there exists a Nil$_s$ Boh$_0$-set $C$
	with the form
	$$C=\bigcap_{k=1}^{t}\{n\in \mathbb{Z}: \{q_k(n)\}\in (-\delta,\delta)\}$$
	such that for all $(i, j)\in \{1, \cdots, d\} \times \{1, \cdots, l\}$, $p_i(n)$ is proper w.r.t. $m_j$ and $C$,	where  $\delta>0$ is a small enough number,
	$s = \mathop {\max }_{1\le i \le d }deg (p_i)$ and	$q_k \in \widehat{SGP_s}, k=1,2,\dots,t$ for some $t\in \mathbb{N}$.
\end{lem}

\begin{rem} \label{rem}
	We call the Nil$_s$ Bohr$_0$-set $C$ above is associated to $\{p_1, \cdots, p_d\}$ and $\{m_1, \cdots, m_l\}$.
\end{rem}

\medskip

\medskip

\section{Proof of Theorem \ref{thm general} for degree $1$ integer-valued generalized polynomials}

In this section, we will prove Theorem \ref{thm general} for degree $1$ integer-valued generalized polynomials. We need  the following lemma.
\begin{lem} \label{lem-N(p,U,V)-deg-1}
	Let $(X, T)$ be a weakly mixing minimal system and
	$p \in \widetilde{SGP_1}$ with $A(p(n))\neq 0$.
	Then for any non-empty open subsets $U, V$ of $X$,
	$$N(p, U, V):=\{n\in \mathbb{Z}: U\cap T^{-p(n)}V \neq \emptyset \}$$
	is thickly syndetic.
\end{lem}
\begin{proof}
	We may assume	
	$p(n)=\sum\limits_{i=1}^{\rm{t_1}} \left\lceil b_i\left\lceil \alpha_i n\right\rceil \right\rceil-\sum\limits_{j=1}^{\rm{t_2}} \left\lceil c_j\left\lceil \beta_j n\right\rceil \right\rceil
	, n\in \mathbb{Z}$ with
	$ t_1, t_2\in \mathbb{N}, \alpha_i, b_i \in \mathbb{R}, i=1, \cdots, t_1$ and $\beta_j, c_j\in \mathbb{R}, j=1, \cdots, t_2$.
	
	Moreover, 
	$$ A(p(n))=\sum_{i=1}^{t_1}b_i \alpha_i- \sum_{j=1}^{t_2}c_j \beta_j \neq 0.$$

	For given non-empty open subsets $U, V$ of $X$,
	since $(X,T)$ is weakly mixing, 
	$$N(U, V):=\{n\in \mathbb{Z}: U\cap T^{-n}V \neq \emptyset \}$$
	is thickly syndetic (see Theorem 4.7 in \cite{HY02} ).
	Then for any $L\in \mathbb{N}$, there exists a syndetic set $A\subset \mathbb{Z}$
	such that
	$$A+ \{0, 1, \cdots, L\} \subset N(U, V).$$
	We denote $A=\{a_1 <a_2< \cdots \}$ and K the gap of $A$.
	Note that for every $n\in \mathbb{Z}$,
	$$\sum_{i=1}^{t_1}b_i(  \alpha_i n-1)-t_1-\sum_{j=1}^{t_2}c_j(  \beta_j n+1)-t_2 <p(n) < \sum_{i=1}^{t_1}b_i(  \alpha_i n+1)+t_1-\sum_{j=1}^{t_2}c_j(  \beta_j n-1)+t_2.$$
	We put $M=\sum_{i=1}^{t_1}b_i\alpha_i-\sum_{j=1}^{t_2}c_j\beta_j, M_0=\sum_{i=1}^{t_1}b_i+\sum_{j=1}^{t_2}c_j+t_1+t_2 $,
	then we have
	$$ Mn-M_0 < p(n)<Mn+M_0.$$
	We can choose $L\in \mathbb{N}$ large enough, such that
	$L \gg  2M_0+8M $.
	
	For $n\in \mathbb{Z}$, if $ p(n) \in \{0, 1,  \cdots, L\}+a_i$ for some $i\in \mathbb{N}$,
	then $U\cap T^{-p(n)}V \neq \emptyset$.
	
	We consider $n\in \mathbb{Z}$ such that
	$$ a_i \le Mn-M_0<  p(n) < Mn+M_0 \le a_i+L$$
	for some $i\in \mathbb{N}$.
	Then we have
	$$ \frac{a_i}{M}+\frac{L}{M}-\frac{M_0}{M}\ge n \ge \frac{a_i}{M}+\frac{M_0}{M} (\ if \ M \ positive),$$
	or
	$$ \frac{a_i}{M}+\frac{L}{M}-\frac{M_0}{M}\le n \le \frac{a_i}{M}+\frac{M_0}{M} (\ if \ M \  negative).$$
	Without loss of generality, we way assume that $M$ is positive.
	
	Since
	$$ \frac{a_i}{M}+\frac{M_0}{M} \le
	\left\lceil\frac{a_i}{M}\right\rceil  +  \left\lceil \frac{M_0}{M}\right\rceil+ 2$$ and
	$$\frac{a_i}{M}+\frac{L}{M}-\frac{M_0}{M} \ge
	\left\lceil\frac{a_i}{M}\right\rceil  + \left\lceil \frac{L}{M}\right\rceil-
	\left\lceil \frac{M_0}{M}\right\rceil -3,$$
	then when
	$$ n\in \{n\in \mathbb{Z}: \left\lceil\frac{a_i}{M}\right\rceil  +  \left\lceil \frac{M_0}{M}\right\rceil+ 2
	\le n\le
	\left\lceil\frac{a_i}{M}\right\rceil  + \left\lceil \frac{L}{M}\right\rceil-
	\left\lceil \frac{M_0}{M}\right\rceil -3 \},$$
	we have that $p(n)\in N(U, V)$.
	
	Let $$B=\{ b_i \buildrel \Delta \over = \left\lceil\frac{a_i}{M}\right\rceil  +  \left\lceil \frac{M_0}{M}\right\rceil+ 2: a_i\in A, i=1,2, \cdots \},$$
	$$L_N={\left\lceil\frac{L}{M}\right\rceil} -
	2\left\lceil\frac{M_0}{M}\right\rceil -5>0. $$
	Then $b_{i+1}-b_i=\left\lceil\frac{a_{i+1}}{M}\right\rceil -\left\lceil\frac{a_i}{M}\right\rceil \le \frac{a_{i+1}}{M}-\frac{a_i}{M}+2
	= \frac{a_{i+1}-a_i}{M}+2\le \frac{K}{M}+2 $ for all $i\in \mathbb{N}$,
	thus $B$ is syndetic.
	Since $L$ can be large enough, so is $L_N.$
	Thus $B+\{0, 1, \cdots, L_N\}\subset N(p, U, V)$, i.e., $N(p, U, V)$ is thickly syndetic.
\end{proof}

\medskip
First we prove an even more special case.

\begin{thm} \label{thm simple}
	Let $(X, T)$ be a weakly mixing minimal system, $p_1, \cdots, p_d\in \widetilde{SGP_1}$ and $(p_1,\cdots,p_d)$ be non-degenerate. 
	Then there is a dense $G_{\delta}$ subset $X_0$ of $X$ such that for all $x\in X_0$,
	$$\{(T^{p_1(n)}x, \cdots, T^{p_d(n)}x): n\in \mathbb{Z}\}$$
	is dense in $X^d$.
	
	Moreover,  for any non-empty open subsets $U, V_1, \cdots, V_d$ of $X$, for any $\varepsilon>0 \ (\varepsilon<\frac{1}{4})$,
	for any $s, t\in \mathbb{N}$ and $g_1, \cdots, g_t \in \widehat{SGP_s}$,
	put
	$$C=C(\varepsilon, g_1, \cdots, g_t)=\bigcap_{j=1}^{t}\{n \in \mathbb{Z}: \{g_i(n)\}\in (-\varepsilon,\varepsilon)\},$$
	$$N=\{n\in \mathbb{Z}: U\cap T^{-p_1(n)}V_1 \cap \cdots \cap T^{-p_d(n)} V_d \neq \emptyset\},$$
	we have  $N \cap C $ is syndetic.

\end{thm}
\begin{proof}	
	By Lemma \ref{lem-equivalent}, it suffices to prove the moreover part of the theorem.
	We will prove it by induction on $d$.
	
	When $d=1$, by Lemma \ref{lem-N(p,U,V)-deg-1}, $N=N(p_1, U,V_1)$ is thickly syndetic, note that $C\in \mathcal{F}_{GP_s}=\mathcal{F}_{s,0}$ is a syndetic set, hence $N\cap C$ is syndetic.

	Assume that the result holds for $d > 1$.
	Next we will show that the result holds for $d+1$.	
	Let $U, V_1, \cdots, V_d, V_{d+1}$ be non-empty open subsets of $X$,
	$ 0<\varepsilon<\frac{1}{4}$ and $g_1, \cdots, g_{t}\in \widehat{SGP_s}$.
	We put	
	$${C}=C(\varepsilon,g_1,\dots,g_t),$$
	$$N=\{n \in \mathbb{Z}: U\cap T^{-p_1(n)}V_1 \cap \dots \cap T^{-p_{d+1}(n)}V_{d+1}\neq \emptyset\},$$
	we will show that $N\cap C$
	is syndetic.
	
	Let	
	$$\tilde{C}_1=C(\frac{\varepsilon}{2}, g_1,\dots,g_t ),$$
	then $\tilde{C}_1 \in \mathcal{F}_{GP_s}=\mathcal{F}_{s,0}$
	is a syndetic set.
	By Lemma \ref{lem-goodm}, there exist $Q\subset \widehat{SGP_b}$ (for some $b\in \mathbb{N}$) and 
	$$\tilde{C}_2=C(\frac{\varepsilon}{2},Q)=\bigcap_{q(n)\in Q}\{n\in \mathbb{Z}:\{q(n)\}\in (-\frac{\varepsilon}{2},\frac{\varepsilon}{2})\}$$
	such that for each $m\in \tilde{C}_2$, $m$ is good w.r.t. 
	$q(n) \in \{p_1(n),p_2(n),\cdots,p_{d+1}(n)\}\cup \{ \left\lceil g_1 \right\rceil, \cdots,\left\lceil g_t \right\rceil \}$. Let $\tilde{C}=\tilde{C}_1\cap \tilde{C}_2$, $\tilde{C}$ is a syndetic set.

	Since $(X, T)$ is minimal, there is some $l\in \mathbb{N}$ such that
	$X=\cup_{j=0}^{l}T^jU$.
	By Lemma \ref{lem-deg-1},
	there are non-empty open subsets $V_1^{(l)}, \cdots, V_{d+1}^{(l)}$
	and integers $k_0, k_1, \cdots, k_l \in  \tilde{C}$
	such that for each $i=1, 2, \cdots, d+1$,
	one has that
	$$T^{p_i(k_j)}T^{-j}V_{i}^{(l)} \subset V_i,  \ \text{for \ all \ }  0\le j \le l.$$

	Notice that $k_i \in \tilde{C} \subset \tilde{C}_2$ is good w.r.t. 	$q(n) \in \{p_1(n),p_2(n),\cdots,p_{d+1}(n)\}\cup \{ \left\lceil g_1 \right\rceil, \cdots,\left\lceil g_t \right\rceil \}$, $0\le i \le l$.
	By Lemma \ref{lem-proper}, there is a Nil$_1$ Bohr$_0$-set $C_1'$ associated to $\{p_1, \cdots, p_{d+1}\}$ and $\{k_0, k_1, \cdots, k_l\}$,
	and by Lemma \ref{asso for special cases}, there is a Nil$_s$ Bohr$_0$-set $C_1''$
	associated to $\{ \left\lceil g_1 \right\rceil, \cdots,\left\lceil g_t \right\rceil \}$ and $\{k_0, k_1, \cdots, k_l\}$.

	Put $C_1=C_1'\cap C_1''$, then $C_1\in \mathcal{F}_{s, 0}$ is a Nil$_s$ Bohr$_0$-set.	We may assume that $\frac{\varepsilon}{2}$ is as in Lemma \ref{asso for special cases}.
	
	Let $q_i=p_{i+1}-p_1 \in \widetilde{SGP_1}$, $i=1, 2, \cdots, d$.
	Then by induction hypothesis,
	$$\{n\in \mathbb{Z}: V_1^{(l)} \cap T^{-q_1(n)}V_2^{(l)}\cap \cdots \cap T^{-q_d(n)}V_{d+1}^{(l)}\neq \emptyset \} \cap (\tilde{C}\cap C_1 )$$
	is syndetic.
	
	Put
	$$E=\{n\in \mathbb{Z}: V_1^{(l)} \cap T^{-q_1(n)}V_2^{(l)}\cap \cdots \cap T^{-q_d(n)}V_{d+1}^{(l)}\neq \emptyset \} \cap (\tilde{C}\cap C_1 ).$$
	Since $E \subset C_1\subset C_1'$,
	we have  	
	$$p_i(m+k_j)= p_i (m)+ p_i(k_j),\forall m \in E $$
	for all $i=1,2,\dots,d+1, j=0,1,\dots,l$.

	Let $m\in E$. Then there is some $x_m \in V_1^{(l)}$ such that
	$T^{q_i(m)}x_m \in V_{i+1}^{(l)}$ for $i=1, \cdots, d$.
	There is some $y_m$ with $y_m=T^{p_1(m)}x_m$.
	Since $X=\cup_{j=0}^{l}T^jU$, there is some $b_m \in \{0, 1, \cdots, l\}$
	such that $T^{b_m}z_m=y_m$ for some $z_m\in U$.
	Thus for each $i=1, 2, \cdots, d+1$,
	\begin{align*}
	T^{ p_i (m+k_{b_m})}z_m
	&=T^{ p_i (m+k_{b_m})}T^{-b_m}y_m \\
	&=T^{ p_i (m+k_{b_m})}T^{-b_m}T^{-p_{1}(m)}x_m \\
	&=T^{ p_i (m)} T^{ p_i (k_{b_m}) }T^{-b_m}T^{-p_{1}(m)}x_m \\
	&=T^{ p_i(k_{b_m})  }T^{-b_m}T^{p_{i}(m)-p_1(m)}x_m \\
	&= T^{p_i(k_{b_m}) }T^{-b_m}T^{q_{i-1}(m)}x_m 	\\
	&\subset T^{p_i (k_{b_m}) }T^{-b_m} V^{(l)}_i  \subset V_i.	\\
	\end{align*}
	That is, 	
	$$z_m \in U\cap T^{-p_1(n)} V_1 \cap \cdots \cap T^{-p_d (n)}V_d  \cap T^{- p_{d+1}(n)} V_{d+1},$$
	where $n=m+k_{b_m} \in N$.
	
	Note that $k_{b_m} \in \tilde{C}$ implies that
	$$\{g_j(k_{b_m})\} \in (-\frac{\varepsilon}{2},\frac{\varepsilon}{2}),$$
	and
	$m \in E \subset C_1''$ implies that
	$$\{g_j(m+k_{b_m})\}\in (   \{g_j (k_{b_m})\}-\frac{\varepsilon}{2}, \{g_j (k_{b_m})\}+\frac{\varepsilon}{2}),$$
	for all $j=1, \cdots, t$.
	Hence $m+k_{b_m} \in C$.
	Thus
	$$N \cap C\supset \{m+k_{b_m}: m\in E\}$$
	is a syndetic set.
	By induction the proof is completed.
\end{proof}
\medskip
Now we can prove our main result for degree $1$ integer-valued generalized polynomials.
\begin{thm}\label{thm special}
	Let $(X, T)$ be a weakly mixing minimal system, $p_1, \cdots, p_d\in \widetilde{GP_1}$ and $(p_1,p_2,\cdots,p_d)$ be non-degenerate.
	Then there is a dense $G_{\delta}$ subset $X_0$ of $X$ such that for all $x\in X_0$,
	$$\{(T^{p_1(n)}x, \cdots, T^{p_d(n)}x): n\in \mathbb{Z}\}$$
	is dense in $X^d$.
	
	Moreover,  for any non-empty open subsets $U, V_1, \cdots, V_d$ of $X$, for any $\varepsilon>0 \ (\varepsilon<\frac{1}{4})$,
	for any $s, t\in \mathbb{N}$ and $g_1, \cdots, g_t \in \widehat{SGP_s}$,
	put
	$$C=C(\varepsilon, g_1, \cdots, g_t)=\bigcap_{j=1}^{t}\{n \in \mathbb{Z}: \{g_i(n)\}\in (-\varepsilon,\varepsilon)\},$$
	$$N=\{n\in \mathbb{Z}: U\cap T^{-p_1(n)}V_1 \cap \cdots \cap T^{-p_d(n)} V_d \neq \emptyset\},$$
	we have  $N \cap C $ is syndetic.
\end{thm}
\begin{proof}
	By Lemma \ref{lem-equivalent}, it suffices to prove the moreover part of the theorem.
	Let $p_1, \cdots, p_d\in \widetilde{GP_1}$. Then by Lemma \ref{lem sgp_gp}, there exists $h_i(n) \in \widetilde{SGP_1}$, $i=1,2,\dots,d$ and  $C_1=C(\delta,q_1,\cdots,q_k)$ such that $p_i(n)=h_i(n), \forall n \in C_1,  i=1,2,\dots,d$.
	
	Set $$N_1=\{n \in \mathbb{Z}: U\cap T^{-h_1(n)}V_1 \cap \cdots \cap T^{-h_d(n)}V_d \neq \emptyset \},$$
	by Theorem \ref{thm simple}, $N_1\cap (C\cap C_1)$ is syndetic. Since for any $n\in N_1\cap (C\cap C_1) \subset C_1$, $p_i(n)=h_i(n),i=1,2,\cdots,d$, we have
	$$N_1 \cap (C \cap C_1) \subset N \cap C,$$
	hence $N\cap C$ is syndetic.
\end{proof}

\medskip
\section{PET-induction and the proof of Theorem \ref{thm general}}

\subsection{The PET-induction}
\ \
In this section, we will prove Theorem \ref{thm general}  using PET-induction,
which was introduced by Bergelson in \cite{Ber-87}.
Basically, we associate any finite collection of integer-valued generalized polynomials with a ``complexity",
and reduce the complexity at some step to the simple one,
where
we use the simple one as the first step (basis of induction).
We first introduce the precise definition of the ``complexity",
in a sense, it is an ordering relationship.

Let $p(n), q(n) \in \widetilde{SGP}$, we denote $p \sim q$ if 
$\deg(p)=\deg(q)$ and $\deg(p-q)<\deg(p)$.
One can easily check that $\sim$ is an equivalence relation. A \emph{system} $P$ is a finite subset of $\widetilde{SGP}$.
Given a system $P$, we define its \emph{weight vector} $\Phi(P)=(\omega_1, \omega_2, \cdots ) $,
where $\omega_i$ is the number of equivalent classes under $\sim$ of degree $i$ integer-valued generalized polynomials represented in $P$.
For distinct weights
$ \Phi(P)= (\omega_1, \omega_2, \cdots)$ and $\Phi(P')= (\upsilon_1, \upsilon_2, \cdots)$,
one writes $\Phi(P)>\Phi(P')$ if $\omega_d>\upsilon_d$, where $d$ is the largest $j$ satisfying
$ \omega_j \neq \upsilon_j$,
then we say that $P'$ \emph{precedes} $P$.
This is a well-ordering of the set of weights and the PET-induction is simply induction on this ordering.

For example, let $P=\{ \left\lceil an\right\rceil +2n, \left\lceil bn^3 \left\lceil cn\right\rceil\right\rceil + \left\lceil en^3  \right\rceil, 4n^4, 4n^4+n^3, \left\lceil fn \right\rceil \left\lceil hn \right\rceil \}$ (where $a, b, c, e, f, h$ are distinct  numbers),
then $\Phi(P)=(1, 1, 0, 2, 0, \cdots )$.    

In order to prove the Theorem \ref{thm general} for system $P$, we will start from $\Phi(P)=(d,0,\cdots)$ (this is true by Theorem \ref{thm simple}). After that, we assume the result holds for any systems $P'$ with $\Phi(P')<\Phi(P)$. Then we show the result holds for $P$, and we complete the proof.

\medskip

\subsection{Some Lemmas}
To simplify the argument, we need to introduce  three symbols: $>>$, $\approx$ and $=_{C}$ .

\begin{itemize}
	\item Let $a>b>0$, we denote $a>>b$ iff there exists a large enough $N>0$  such that $a>N(b+1)$.
	\item $a \approx b$ iff $|a|>>|a-b|$ and $|b|>>|a-b|$.	
	\item $p(n)=_{C}q(n)$ iff $p(n)=q(n)$ for any $n \in C$.
\end{itemize}

We have the following observation.
\begin{lem} \label{lem-large-approx}
	\begin{enumerate} 
		\item Let $|a|>> 1$. Then $\left\lceil a \right \rceil \approx a$.
		\item Let $|a|>>1$, $|b|>>1$, $a\approx a'$ and $b \approx b'$. Then
		$ab \approx a'b'$.
		Moreover, if it still satisfies $|a+b|>>1$, then $a'+b'\approx a+b$.
		\item Let $|\sum_{i=1}^{k}a_i|>>1$, and for any $1 \le i \le k$, $|a_i|>>1$, $a_i\approx a_i'$.  Then $|\sum_{i=1}^ka_i'|>>1$.
		
	\end{enumerate}
\end{lem}

\medskip

For instance,
$10000\sqrt{2}>>1$, $5000\sqrt{3}>>1$, $|10000\sqrt{2}-5000\sqrt{3}|>>1$, 
$\left \lceil 10000\sqrt{2} \right\rceil \approx 10000\sqrt{2}, 
\left \lceil 5000\sqrt{3} \right\rceil \approx 5000\sqrt{3},$
$10000\sqrt{2} \times 5000\sqrt{3} \approx \left \lceil 10000\sqrt{2} \right\rceil \left \lceil 5000\sqrt{3} \right\rceil.$

\medskip

Recall that $A(p(n))$ be the sum of the coefficients of the maximal-degree components of the generalized polynomial $p(n)$ (see Definition \ref{def-nondegenerate}).  We have the following lemmas.

\begin{lem}\label{lem-A(p)-simple}
	If $h(n)=\left \lceil p(n) \right \rceil\in \widetilde{SGP_d}$ with $\deg(h)\ge 2$ and 
	$p(n)=\prod_{i=1}^{k}L(a_{i,1}n^{j_{i,1}},\cdots, a_{i,l_{i}}n^{j_{i,l_i}})$
	where $|a_{i,1}|>>1$, $j_{i,1} \ge 0$, $a_{i,1}\in \mathbb{R}$ and
	$|a_{i,t}|>>1$, $j_{i,t} \ge 1$, $a_{i,t} \in \mathbb{R}\setminus \mathbb{Q}$ for $t=2,\cdots,l_i$, $1 \le i \le k$.
	Let $0\neq m \in \mathbb{Z}$. Then there exist a Nil Bohr$_0$ set $C$ and 
	$D(h(n),m) \in \widetilde{SGP_{d-1}}$ such that
	$$D(h(n),m)=_{C}h(n+m)-h(n)-h(m)$$ and
	$$A(D(h(n),m))\approx \deg(h) mA(h(n)).$$
	We call $D(h(n),m)$ be the deritive of $h(n)$ w.r.t. $m$.
	
\end{lem}
\begin{proof}
	We first prove it for two special cases to illustrate the idea. Then we prove it for the general case.
	
	Note: For any $k\in \mathbb{N}$, $k\mathbb{Z}=\{n\in \mathbb{Z}: \{\frac{n}{k} \} \in (-\frac{1}{2k}, \frac{1}{2k} )\}$ is a Nil$_1$ Bohr$_0$-set, and for any $a\in \mathbb{Q}$, there exists $k_0\in \mathbb{N}$ such that $\left\lceil an\right\rceil=_{k_0\mathbb{Z}} an $. 
	
	\noindent{\bf Special case 1.} Assume that $p(n)=L(an,bn^2)$ with $a>>1$, $b>>1$ and $a\in \mathbb{R}, b \in \mathbb{R}\setminus \mathbb{Q}$.
	By Note, we may assume that $a \in \mathbb{R}\setminus \mathbb{Q}$.
	Let $0\neq m \in \mathbb{Z}$, then $m$ is good w.r.t. $p(n)$. The expansion of
	$b(n+m)^2$ is
	$$b(n+m)^2=b\sum_{i=0}^2C_2^im^in^{2-i}=b(n^2+2mn+m^2).$$
	By Lemma \ref{lem-proper-simple}, there exists a Nil Bohr$_0$ set $C$ such that
	$$\left \lceil b(n+m)^2 \right\rceil =_{C} \left \lceil bn^2\right\rceil + \left \lceil 2bmn\right\rceil +\left \lceil bm^2\right\rceil,$$
	\begin{eqnarray*}
		\left \lceil a(n+m) \left \lceil b(n+m)^2 \right\rceil \right\rceil
		&=_{C}& \left \lceil (an+am) (\left \lceil bn^2\right\rceil + \left \lceil 2bmn\right\rceil +\left \lceil bm^2\right\rceil )\right\rceil \\
		&=_{C}&\left\lceil an \left \lceil bn^2\right\rceil \right\rceil  
		+ \left\lceil an \left\lceil 2bmn\right\rceil \right\rceil + \left\lceil an\left \lceil bm^2\right\rceil \right\rceil \\
		&&+\left\lceil am \left \lceil bn^2\right\rceil \right\rceil  
		+ \left\lceil am \left \lceil 2bmn \right\rceil \right\rceil 
		+ \left\lceil am\left \lceil bm^2\right\rceil \right\rceil 
	\end{eqnarray*}
	We denote
	\begin{eqnarray*}
		D(h(n),m)&=&\left\lceil an \left\lceil 2bmn\right\rceil \right\rceil + \left\lceil an\left \lceil bm^2\right\rceil \right\rceil+ 
		\left\lceil am \left \lceil bn^2\right\rceil \right\rceil  
		+ \left\lceil am \left \lceil 2bmn \right\rceil \right\rceil, 
	\end{eqnarray*}
	then $D(h(n),m)=_{C}h(n+m)-h(n)-h(m)$. The maximal-degree components of $D(h(n),m)$ are $\left\lceil an\left \lceil 2bmn\right\rceil \right\rceil$ and
	$\left\lceil am \left \lceil bn^2\right\rceil \right\rceil$, hence
	$$A(D(h(n),m))=2abm+abm=3abm=\deg(h)mA(h(n)).$$
	
	\noindent{\bf Special case 2.} Assume that $p(n)=\left \lceil an \right\rceil \left \lceil bn \right\rceil $.

	By Note, we may assume that $a, b\in \mathbb{R}\backslash \mathbb{Q}$.
	Let $0\neq m\in \mathbb{Z}$, then $m$ is good w.r.t. $p(n)$. By Lemma \ref{lem-proper-simple}
	there exist a Nil Bohr$_0$ set $C$ such that
	\begin{eqnarray*}
		\left \lceil a(n+m) \right\rceil \left \lceil b(n+m)\right\rceil
		&=_{C}& (\left \lceil an \right\rceil+ \left \lceil am \right\rceil) 
		(\left \lceil bn \right\rceil +\left \lceil bm \right\rceil) \\
		&=& \left \lceil an \right\rceil \left \lceil bn \right\rceil+ \left \lceil am \right\rceil \left \lceil bn \right\rceil
		+\left \lceil an \right\rceil\left \lceil bm \right\rceil +\left \lceil am \right\rceil\left \lceil bm \right\rceil.
	\end{eqnarray*}
	Let $D(h(n),m)=\left \lceil am \right\rceil \left \lceil bn \right\rceil
	+\left \lceil an \right\rceil\left \lceil bm \right\rceil$, then $D(h(n),m)=_{C}h(n+m)-h(n)-h(m)$ and 
	$$A(D(h(n),m))=a\left \lceil bm \right\rceil +\left \lceil am \right\rceil b \approx 2mab=\deg(h)mA(h(n)).$$
	
	\noindent{\bf The general case.} Assume that $p(n)=\prod_{i=1}^{k}L(a_{i,1}n^{j_{i,1}},\cdots, a_{i,l_{i}}n^{j_{i,l_i}})$.
	By the argument of Special case 1 and 2, we may assume that for $0\neq m \in \mathbb{Z}$, $m$ is good w.r.t. $p(n)$.
	
	By Lemma \ref{lem-proper-simple}, there exist a Nil Bohr$_0$ set $C$ and
	$D(h(n),m)=_{C}h(n+m)-h(n)-h(m)$. The maximal-degree components of $D(h(n),m)$
	are 
	$$\left \lceil \prod_{i=1}^{k}L(a_{i,1}n^{j_{i,1}},\cdots,a_{i,t}C_{j_{i,t}}^{1}mn^{j_{i,t}-1},\cdots, a_{i,l_{i}}n^{j_{i,l_i}})\right\rceil,  1 \le t \le l_i, 1\le i \le k.$$
	Hence
	\begin{eqnarray*}
		A(D(h(n),m)) 
		&\approx& \sum_{i=1}^{k}\sum_{t=1}^{l_i}(a_{i,1}\cdots (C_{j_{i, t}}^1ma_{i,t})\cdots a_{i,l_i})\prod_{s\neq i,s=1}^{k}(a_{s,1}\cdots  a_{s,l_s})\\
		&=&\prod_{s=1}^{k}(a_{s,1}\cdots  a_{s,l_s})\sum_{i=1}^{k}\sum_{t=1}^{l_i} j_{i, t}m \\
		&=& A(h(n))\deg(h)m.
	\end{eqnarray*}
\end{proof}

\medskip

\begin{lem} \label{lem-A(p)-general}
	Let $h(n) =\sum_{k=1}^l c_k\left \lceil p_k(n) \right \rceil \in \widetilde{SGP_d}$, where $c_k \in \mathbb{Z}$, $p_{k}(n) \in \widehat{SGP_d}$  as in above lemma, $|A(h(n))|>>1$ and $\deg(h)\ge 2$. Let $0\neq m\in \mathbb{Z}$. Then there exist a Nil Bohr$_0$ set $C$ and 
	$D(h(n),m) \in \widetilde{SGP_{d-1}}$ such that
	$$D(h(n),m)=_{C}h(n+m)-h(n)-h(m),$$ 
	$$A(D(h(n),m))\approx \deg(h) mA(h(n)).$$
	
\end{lem}
\begin{proof}
	Notice that when calculating $A(h)$, we just need to consider the maximal-degree components, we can assume that $\deg(p_k(n))=\deg(h(n)), k=1, \cdots, l$. Then $A(h(n))=\sum_{k=1}^{l}A(p_k(n))$.
	For any $k=1, \cdots, l$,
	by Lemma \ref{lem-A(p)-simple}, there exist Nil Bohr$_0$ set $C_k$ 
	and $D(\left \lceil p_k(n)\right\rceil,m)$ such that
	$$D(\left \lceil p_k(n)\right\rceil,m)=_{C_k}\left \lceil p_k(n+m)\right\rceil-
	\left \lceil p_k(n)\right\rceil-\left \lceil p_k(m)\right\rceil,$$
	$$A(D(\left \lceil p_k(n)\right\rceil,m))\approx \deg(p_k)mA(p_k(n)).$$
	Let $C =\bigcap_{k=1}^l C_k$ and 
	$D(h(n),m)=\sum_{k=1}^lc_k D(\left \lceil p_k(n)\right\rceil,m)$. Then 
	$$D(h(n),m)=_{C}h(n+m)-h(n)-h(m),$$
	$$A(D(h(n),m))=\sum_{k=1}^lc_k A(D(\left \lceil p_k(n)\right\rceil,m))\approx \sum_{k=1}^lc_k \deg(p_k)mA(p_k(n))=\deg(h)m A(h(n)).$$
\end{proof}

\begin{lem}\label{lem-h_1-h_2}
	Let $h_1,h_2 \in \widetilde{SGP}$, $(h_1,h_2)$ be non-degnerate, $h_1 \sim h_2$, $\deg(h_1) \ge 2$ and $h_1,h_2 $ satisify conditions in the above lemmas. Then for any $0\neq m \in \mathbb{Z}$,
	$$D(h_1(n),m)-D(h_2(n),m)=D(h_1(n)-h_2(n),m).$$
\end{lem}	
\begin{proof}
	Let $\deg(h_1)=d$. 
	Since $h_1 \sim h_2$, then $\deg (h_1-h_2)=r<d$.
	Then if we write
	$$h_1(n)=\sum_{k=1}^{d}\sum_{j=1}^{l_{k,1}}\left \lceil p_{k,j,1}(n)\right\rceil,\ \  h_2(n)=\sum_{k=1}^{d}\sum_{j=1}^{l_{k,2}} \left \lceil p_{k,j,2}(n) \right\rceil$$
	with $p_{k,j,i} \in \widehat{SGP_k}$, $\deg(\left \lceil p_{k,j,i} \right\rceil ) =k, k=1,2,\dots,d, j=1,2,\cdots l_{k,i},  i=1,2.$
	Then for each $r+1 \le k \le d$, 
	$$\sum_{j=1}^{l_{k,1}}\left \lceil p_{k,j,1}(n)\right\rceil=\sum_{j=1}^{l_{k,2}}\left \lceil p_{k,j,2}(n)\right\rceil.$$
	Hence
	$$D(h_1(n),m)-D(h_2(n),m)=D(h_1(n)-h_2(n),m).$$
\end{proof}

\begin{lem}\label{lem-multi-q(ij)}
	Let $p_i \in \widetilde{SGP}$, $(p_1,p_2,\cdots,p_d)$ be non-degenerate with  $\deg(p_i)\ge 2, 1\le i \le d$ and $p_i$ satisify conditions in the above lemmas. Then there exist a
	sequence $\{r(n)\}_{n=0}^{\infty}$ of natural numbers, such that for any $l \in \mathbb{N}$ and $k_0,k_1,\cdots,k_l \in \mathbb{Z}$	with $|k_0|>r(0)$ and $|k_{i}|>|k_{i-1}|+r(|k_{i-1}|)$,
	there exist a Nil Bohr$_0$ set $C$ and $q_{i,j}\in \widetilde{SGP}$ with
	$$q_{i,j}(n):=D(p_i(n),k_j)+p_i(n)-p_1(n)=_Cp_i(n+k_j)-p_i(k_j)-p_1(n)$$
	and
	\begin{equation}\label{eq-a(qij-q)}
	|A(q_{i,j}(n))|>>1, |A(q_{i,j}(n)-q_{i',j'}(n))|>>1 
	\end{equation}
	for all $ (i,j) \neq (i',j')\in \{1,2,\cdots,d\}\times\{0,1,\cdots,l\}$.
\end{lem}
\begin{proof}	
	Let $$M=\max_{1 \le i \neq i' \le d}\{\deg(p_i),|A(p_i)|, |A(p_i-p_{i'})|\},$$
	$$L=\min_{1 \le i \neq i' \le d}\{\deg(p_i),|A(p_i)|, |A(p_i-p_{i'})|\}$$
	Set $r(n)=10^{10}\frac{M^2}{L^2}(n+1), n=0, 1, \cdots$. We will show that if
	$|k_i|>|k_{i-1}|+r(|k_{i-1}|)$, then for all $ (i,j) \neq (i',j')\in \{1,2,\cdots,d\}\times\{0,1,\cdots,l\}$, \eqref{eq-a(qij-q)} holds.
	To do so, for $k_j,k_{j'} \in \mathbb{Z}$, we need to calculate
	$A(q_{i,j}(n))$ and $A(q_{i,j}(n)-q_{i',j'}(n))$.
	
	\noindent{\bf Case 1: The value of  $A(q_{i,j}(n))$.} 
	Notice that $q_{i,j}(n)=D(p_i(n), k_j)+p_i(n)-p_1(n)$. 
	\begin{itemize}
		\item If $p_i(n)\nsim p_1(n)$,
		then the maximal-degree components of $q_{i,j}(n)$ is either in $p_i(n)$ or in $p_1(n)$,
		hence $A(q_{i,j}(n))$ is equal to $A(p_i(n))$ or 
		$A(p_1(n))$, hence $|A(q_{i,j}(n))|>>1$. 
		\item 	If $p_i(n) \sim p_1(n)$, there are two cases. 
		If $\deg(p_i-p_1)<\deg(p_i)-1$, then
		\begin{equation}\label{eq-q-1}
		A(q_{i,j}(n))=A(D(p_i,k_j))\approx k_j\deg(p_i)A(p_i(n)). 
		\end{equation}
		If $\deg(p_i-p_1)=\deg(p_i)-1$, then if $|k_j\deg(p_i)A(p_i(n))+A(p_i(n)-p_1(n)) |>>1$ and by Lemma \ref{lem-large-approx}, we have
		\begin{equation}\label{eq-q-2}
		A(q_{i,j}(n))\approx k_j\deg(p_i)A(p_i(n))+A(p_i(n)-p_1(n)).
		\end{equation} 
	\end{itemize}
	
	\noindent{\bf Case 2: The value of $A(q_{i,j}(n)-q_{i', j'}(n))$.}
	$$q_{i,j}(n)-q_{i',j'}(n)=D(p_i,k_j)-D(p_{i'},k_{j'})+p_{i}(n)-p_{i'}(n).$$
	\begin{itemize}
		\item  If $p_i(n) \nsim p_{i'}(n)$, then $A(q_{i,j}(n)-q_{i',j'}(n))$ is equal to $A(p_i(n))$ or $A(p_{i'}(n))$, hence $|A(q_{i,j}(n))|>>1$ for all $j=0, \cdots, l$.
		\item  If $p_i(n) \sim p_{i'}(n)$, there are two cases. If $j=j'$, then by Lemma \ref{lem-h_1-h_2},
		\begin{equation*} 
		D(p_{i}(n),k_j)-D(p_{i'}(n),k_j)= D(p_{i}(n)-p_{i'}(n),k_j)
		\end{equation*}
		and hence $|A(q_{i,j}(n)-q_{i',j}(n))|=|A(p_{i}(n)-p_{i'}(n))|>>1.$
		If $j \neq j'$, there are two cases.  
		If $\deg(p_i-p_{i'})<\deg(p_i)-1$, then if $|k_j\deg(p_i)A(p_i)-k_{j'}\deg(p_{i'})A(p_{i'}) |>>1$ and by Lemma \ref{lem-large-approx}, one has 
		\begin{equation}\label{eq-qq1}
		A(q_{i,j}(n)-q_{i',j'}(n))\approx k_j\deg(p_i)A(p_i)-k_{j'}\deg(p_{i'})A(p_{i'}).
		\end{equation}
		If $\deg (p_i-p_{i'})=\deg (p_i)-1$, then if $|k_j\deg(p_i)A(p_i)-k_{j'}\deg(p_{i'})A(p_{i'})+A(p_i-p_{i'}) |>>1$ and by Lemma \ref{lem-large-approx}, one has 
		\begin{equation}\label{eq-qq2}
		A(q_{i,j}(n)-q_{i',j'}(n)) \approx k_j\deg(p_i)A(p_i)-k_{j'}\deg(p_{i'})A(p_{i'})+A(p_i-p_{i'}).
		\end{equation}
	\end{itemize}
	
	\medskip
	Now we will show that \eqref{eq-a(qij-q)} holds. 
	First choose any $|k_0|>r(0)$, then by \eqref{eq-q-1}  $$|A(q_{i,0})|>>1,$$ by \eqref{eq-q-2}
	$$|A(q_{i,0})|\ge |k_0|L^2-M>10^{10}\frac{M^2}{L^2}L^2-M>>1,$$
	Thus  \eqref{eq-a(qij-q)} holds for $k_0$.
	
	Next we choose $|k_1|>|k_0|+r(|k_0|)$, by \eqref{eq-q-1} and 
	\eqref{eq-q-2},
	$$|A(q_{i,0})|>>1, |A(q_{i,1})|>>1.$$
	By \eqref{eq-qq1} $$|A(q_{i,1}-q_{i',0})|\ge |k_1| L^2-|k_0|M^2>(|k_0|+10^{10}\frac{M^2}{L^2}(|k_0|+1))L^2-|k_0|M^2>>1.$$
	By \eqref{eq-qq2}
	$$|A(q_{i,1}-q_{i',0})|>|k_1|L^2-|k_0|M^2-M>(|k_0|+10^{10}\frac{M^2}{L^2}(|k_0|+1))L^2-|k_0|M^2-M>>1.$$
	Hence \eqref{eq-a(qij-q)} holds for $k_0,k_1$.
	
	Inductively, we can show that \eqref{eq-a(qij-q)} holds for $k_0,k_1,\cdots,k_l$, and we complete the proof.
	
\end{proof}

\medskip

\begin{lem}\label{lem-g(n)-g(kn)}
	Let $U, V_1, \cdots, V_d$ be non-empty open sets of $X$. Let $k \in \mathbb{N}$ and we denote
	$$N=\{n\in \mathbb{Z}: U\cap T^{-p_1(n)}V_1 \cap \cdots \cap T^{-p_d(n)} V_d \neq \emptyset\},$$
	$$N_1=\{m\in \mathbb{Z}: U\cap T^{-p_1(km)}V_1 \cap \cdots \cap T^{-p_d(km)} V_d \neq \emptyset\}.$$
	If for any $\varepsilon>0$, for any $s, t \in \mathbb{N}$ and $g_1, \cdots, g_t\in \widehat{SGP_s} $, 
	$$N_1\cap C(\varepsilon, g_1, \cdots, g_t )$$ is syndetic, 
	then $N\cap  C(\varepsilon, g_1, \cdots, g_t )$ is syndetic.	
\end{lem}
\begin{proof}
	Let $\tilde{g_i}(n)=g_i(kn) \in \widehat{SGP_s}$ and $$C_1=C(\varepsilon,\tilde{g}_1,\cdots,\tilde{g}_t)\cap C(\varepsilon,g_1,\cdots,g_t),$$
	$N_1\cap C_1$ is syndetic. For any $n \in N_1 \cap C_1$, 
	$kn \in N\cap C(\varepsilon,g_1,\cdots,g_t)$.
	Since  $N_1 \cap C_1$ is syndetic, $N\cap C(\varepsilon,g_1,\cdots,g_t)$ is syndetic.
\end{proof}

\medskip
\subsection{The proof of Theorem \ref{thm general}}
Notice that for any $0\neq a \in \mathbb{R}$, if we choose $$C=\{n \in \mathbb{Z}: \{abn^k\}\in (-\frac{1}{4}, \frac{1}{4}), \{bn^k\}\in (-\frac{1}{4|a|},-\frac{1}{4|a|})\},$$ 
then we have $\left\lceil a \left \lceil bn^k\right \rceil\right\rceil=_C \left \lceil abn^k\right\rceil$. Combining this fact with Lemma \ref{lem-g(n)-g(kn)}, from now on we always assume that all $p(n) \in \widetilde{SGP}$ in the following theorem satisfy the conditions in Lemma \ref{lem-A(p)-simple} and Lemma \ref{lem-A(p)-general}.

We first prove the following theorem.

\begin{thm} \label{thm general sgp}
	Let $(X, T)$ be a weakly mixing minimal system, $p_1, \cdots, p_d\in \widetilde{SGP}$ and $(p_1,p_2,\cdots,p_d)$ be non-degenerate.
	Then there is a dense $G_{\delta}$ subset $X_0$ of $X$ such that for all $x\in X_0$,
	$$\{(T^{p_1(n)}x, \cdots, T^{p_d(n)}x): n\in \mathbb{Z}\}$$
	is dense in $X^d$.
	
	{Moreover, for any non-empty open subsets $U, V_1, \cdots, V_d$ of $X$, for any $\varepsilon>0$, for any $s, t \in \mathbb{N}$ and $g_1, \cdots, g_t\in \widehat{SGP_s} $, let
		$$C=C(\varepsilon, g_1, \cdots, g_t ),$$
		$$N=\{n\in \mathbb{Z}: U\cap T^{-p_1(n)}V_1 \cap \cdots \cap T^{-p_d(n)} V_d \neq \emptyset\},$$
		we have $N \cap C $ is syndetic.}
\end{thm}

\begin{proof}
	We will use the PET-induction.  Let $P=\{p_1, \cdots, p_d\}$. Just as the argument above (by Lemma \ref{lem-g(n)-g(kn)}), we can assume that
	$|A(p_i)|>>1$ and $|A(p_i-p_j)|>>1$ for any $1 \le i\neq j \le d$.

	We start from the system whose weight vector is $(d, 0, \cdots)$.
	That is, the degree of all the elements of $P$ is 1.
	By Lemma \ref{lem-N(p,U,V)-deg-1} and Theorem \ref{thm simple}, we know that
	
	\begin{itemize}
		\item[$*_1$] $(X, T)$ is $P$-thickly-syndetic transitive.
		\item[$*_2$] For any non-empty open subsets $U, V_1, \cdots, V_d$ of $X$,
		for any $\varepsilon>0$, for any $s, t\in \mathbb{N}$ and
		$g_1, g_2, \cdots, g_t \in \widehat{SGP_s}$,
		put
		$$C=C(\varepsilon, g_1, \cdots, g_t ),$$
		$$N=\{n\in \mathbb{Z}: U\cap T^{-p_1(n)}V_1 \cap \cdots \cap T^{-p_d(n)} V_d \neq \emptyset\},$$
		we have $N \cap C $ is syndetic.
	\end{itemize}
	
	Let $P \subset \widetilde{SGP}$ be a system whose weight vector $> (d, 0, \cdots)$,
	and we assume that for all systems $P'$ preceding $P$ satisfy $*_1$ and $*_2 $.
	
	Now we show that system $P$ holds. More precisely, 
	in Claim 1 we will show that $*_1$ holds for $P'$ and $*_2$ hold for $P'$ imply that $*_1$ holds for $P$, in Claim 2 we will show that $*_1$ holds for $P$ and $*_2$ holds for $P'$ imply that $*_2$ holds
	for $P$.
	\medskip
	
	{\noindent\bf Claim 1.} \label{ main claim 1}
	$*_1$ holds for $P$, i.e. $(X, T)$ is $P$-thickly-syndetic transitive.
	
	{\noindent\bf Proof of Claim 1:}
	Since the intersection of two thickly syndetic sets is still a thickly syndetic set, it is sufficient to show that for any $p\in P$, and for any given non-empty open subsets $U, V$ of $X$,
	$$N(p,U,V)=\{n\in \mathbb{Z}: U\cap T^{-p(n)}V \neq \emptyset\}$$
	is thickly syndetic.
	
	If $\deg(p)=1$, by Lemma \ref{lem-N(p,U,V)-deg-1}, $N(p,U,V)$ is thickly syndetic.
	
	We assume $\deg(p) \ge 2$.	
	As $(X, T)$ is minimal, there is some $l\in \mathbb{N}$ such that $X=\cup_{i=0}^{l}T^iU$.
	
	Let $L\in \mathbb{N}$ and $k_i=i(L+2)+1$, for all $i\in \{0, 1, \cdots, l\}$.
	Since $(X, T)$ is weakly mixing and minimal, for any  $(i,j) \in \left\{ {0,1, \cdots ,l} \right\}\times \{ {0,1, \cdots ,L} \}$,  $N(V,(T^{p(k_i+j)-i})^{-1}V)$ is thickly syndetic (see Theorem 4.7 in \cite{HY02} ), hence
	$$C:= \bigcap\limits_{(i,j) \in \left\{ {0,1, \cdots ,l} \right\} \times \left\{ {0,1, \cdots ,L} \right\}} \{ k\in \mathbb{Z}: V\cap T^{-k}(T^{p(k_i+j)-i})^{-1}V \neq \emptyset\}$$
	is a thickly syndetic set.
	Choose $c\in C$.
	Then for any $(i,j) \in \{ {0,1, \cdots ,l} \} \times \{ {0,1, \cdots ,L} \} $,
	one has
	$$ V_{i, j}:= V\cap (T^{p(k_i+j)+c-i})^{-1}V $$
	is a non-empty open subset of $V$ and
	$$T^{p(k_i+j)+c-i} V_{i, j} \subset V.  $$
	By Lemma \ref{lem-A(p)-general} and Lemma \ref{lem-proper},
	there is a Nil$_{\deg(p)}$ Bohr$_0$-set $C_1$ associated to $p$ and $\{k_i+j: 0\le i \le l , 0\le j \le L \}$.
	This means for every $(i,j) \in \{ {0,1, \cdots ,l} \} \times \{ {0,1, \cdots ,L} \}$,
	there exists $D(p(n),k_i+j)=_{C_1} p(k_i+j+n)-p(k_i+j)-p(n)$ with $\deg(D(p(n),k_i+j)) <\deg(p)$. Let $q_{i,j}(n)=D(p(n),k_i+j)$ and
	$$P'=\{q_{i, j}: (i,j) \in \{0,1, \cdots ,l\} \times \{ 0,1, \cdots ,L\} \}.$$
	Then $P' \subset \widetilde{SGP}$. 
	Since for any $q_{i, j}\in P'$, $\deg(q_{i,j})<\deg(p)$, we have $\Phi(P') < \Phi(\{p\})$.
	
	For any $(i,j)\neq (i',j')\in \{ {0,1, \cdots ,l} \} \times \{ {0,1, \cdots ,L} \} $, recall that we choose $k_i=i(L+2)+1$,
	$k_i+j \neq k_{i'}+j'$. Hence by Lemma \ref{lem-A(p)-general} and Lemma \ref{lem-large-approx},
	$$A(q_{i,j})\approx\deg(p)(k_i+j)A(p(n)), \ |A(q_{i,j})|>>1,$$ 	$$A(q_{i',j'})\approx\deg(p)(k_{i'}+j')A(p(n)),\  |A(q_{i',j'})|>>1,$$
	$$A(q_{i',j'}-q_{i,j})\approx\deg(p)(k_{i'}+j'-k_i-j)A(p(n)),\ |A(q_{i',j'}-q_{i,j})|>>1.$$

	By the inductive assumption that  $*_2$ holds for $A'$,
	we have
	$$E=\{n\in \mathbb{Z}: V\cap \bigcap_{(i,j) \in \{ {0,1, \cdots ,l} \} \times \{ {0,1, \cdots ,L} \}}T^{-q_{i, j}( n)}V_{i, j}  \neq \emptyset\} \cap C_1$$
	is syndetic.
	
	For $m\in E$, we have $q_{i, j}(m)=p(k_i+j+m)-p(k_i+j)-p(m)$.
	And there exists $x_m\in V$ such that $T^{q_{i,j}(m)}x_m \in V_{i, j}$ for all
	$(i,j) \in \{ {0,1, \cdots ,l} \} \times \{ {0,1, \cdots ,L} \} $.
	Let $y_m=T^{-p(m)}x_m$. Since $X=\cup_{i=0}^{l}T^iU$,
	there are $z_m\in U$ and $0\leq b_m \leq l$
	such that $T^cy_m=T^{b_m}z_m$.
	Then $z_m=T^{-p(m)+c-b_m}x_m$ and we have
	\begin{align*}
	T^{p(m+k_{b_m}+j)}z_m
	&=T^{p(m+k_{b_m}+j)}T^{-p(m)+c-b_m}x_m\\
	&=T^{p(k_{b_m}+j)+c-b_m}(T^{p(m+k_{b_m}+j)-p(k_{b_m}+j)-p(m)}x_m)\\
	&=T^{p(k_{b_m}+j)+c-b_m}(T^{q_{b_m,j}(m)}x_m)\\
	& \in T^{p(k_{b_m}+j)+c-b_m}V_{b_m,j} \subset V
	\end{align*}
	for each $j\in \{0, 1, \cdots, L\}$.
	Thus
	$$\{m+k_{b_m}+j: 0\leq j\leq L\}\subset N(p, U, V).$$
	Hence the set
	$\{n\in \mathbb{Z}: n+j\in N(p, U, V) \ \text{for}\ j=0, 1, \cdots, L \}$
	contains the syndetic set $\{m+k_{b_m}: m\in E\}$.
	As $L\in \mathbb{N}$ can be arbitrary large, $N(p, U, V)$ is a thickly syndetic set.

	\medskip
	
	{\noindent\bf Claim 2.}	$*_2$ holds for $P$. That is, for any non-empty open subsets $U, V_1, \cdots, V_d$ of $X$,
	for any $\varepsilon>0 $, for any $s, t\in \mathbb{N}$  and
	$g_1, g_2, \cdots, g_t \in \widehat{SGP_s}$ ,
	put
	$$C=C(\varepsilon, g_1, \cdots, g_t),$$
	$$N=\{n\in \mathbb{Z}: U\cap T^{-p_1(n)}V_1 \cap \cdots \cap T^{-p_d(n)} V_d \neq \emptyset\},$$
	we have $N \cap C $ is syndetic.
	
	{\noindent\bf Proof of Claim 2:}
	By permuting the indices, we may assume that $\deg(p_i)$ will
	not decrease as $i$ increase. Assume that $\deg(p_w)=1$ and $\deg(p_{w+1})\ge 2$, $1\le w<d$. 
	If for any $p\in P$, $\deg p\ge 2$, we put $w=0$.
	Let $\{r(n)\}_{n=0}^{\infty}$ be the sequence in Lemma \ref{lem-multi-q(ij)} w.r.t. $(p_{w+1},\cdots,p_d)$.

	Put
	$$\tilde{C}=C(\frac{\varepsilon}{2}, g_1, \cdots, g_t),$$
	$$h_1=\max_{p\in A} \text{deg} p, \ h_2= \max_{1\le j\le t} \text{deg} g_j.$$

	Since $(X, T)$ is minimal, there is some $l\in \mathbb{N}$ such that $X=\cup_{i=0}^{l}T^iU$.

	By Claim 1, $(X,T)$ is P-thickly-syndetic transitive.
	Then by Lemma \ref{lem-deg-1},
	there are integers $\{k_j\}_{j=0}^{l} \subset \tilde{C}$
	and non-empty open sets $V_i^{(l)}\subset V_i, 1\leq i\leq d$
	such that
	$|k_j| >|k_{j-1}|+r(|k_{j-1}|)$ for $j=0, \cdots, l \ (k_{-1}=0)$
	and
	$$ T^{p_i(k_j)}T^{-j}V_i^{(l)} \subset V_i, 0\leq j\leq l, 1\leq i\leq d. $$
	By Lemma \ref{lem-proper}, there is a Nil$_{h_1}$ Bohr$_0$-set $C_1'$  associated to $\{p_1, \cdots, p_d\}$ and $\{k_0, \cdots, k_l\}$.
	By Lemma \ref{asso for special cases}, there is a Nil$_{h_2}$ Bohr$_0$-set $C_1''$  associated to $\{g_1, \cdots, g_t\}$ and $\{k_0, \cdots, k_l\}$.
	Put $C_1=C_1'\cap C_1''$, then $C_1\in \mathcal{F}_{h, 0}$, where $h=\max \{h_1, h_2\}$.
	Without loss of generality, we may assume that $\frac{\varepsilon}{2}$ is as in Lemma \ref{asso for special cases}.
	
	Fix $(i, j)\in \{1, \cdots, d\} \times \{0, \cdots, l\}$. 
	For $w+1\le i\le d$, by Lemma \ref{lem-A(p)-general}, 
	there exists $D(p_i(n),k_j) \in \widetilde{SGP}$ with $\deg(D(p_i(n),k_j))<\deg(p_i)$ such that
	$$D(p_i(n),k_j) =p_i(k_j+n)-p_i(k_j)-p_i(n), \forall n\in C_1.$$
	Let $p_{i, j}(n)=p_i(k_j+n)-p_i(k_j)-p_1(n)$ and
	$q_{i, j}(n)=D(p_i(n),k_j)+p_i(n)-p_1(n)$, then $q_{i,j}(n) \in \widetilde{SGP}$ and $$p_{i,j}(n)=q_{i,j}(n), \forall n\in C_1.$$

	For $w\ge 1, 1 \le i \le w$, since $\deg(p_i)=1$,
	we have $p_i(k_j+n)-p_i(k_j)-p_i(n)=0, \forall n\in C_1$,  
	and we let $q_{i,j}(n)=p_i(n)-p_1(n), j=0,1,\cdots,l$.
	
	Let $$P'=\{p_2(n)-p_1(n),\cdots,p_w(n)-p_1(n)\}\cup \{q_{i,j}(n): (i,j)\in \{w+1,\cdots, d\}\times\{0,1,\cdots,l\}\}.$$
	Then $P'\subset \widetilde{SGP}$ and $\Phi(P')<\Phi(P)$ since $q_{i, j}\sim p_i, (i,j)\in \{w+1,\cdots, d\}\times\{0,1,\cdots,l\}$.
	
	For $w=0$, one has $q_{1, j}(n)=D(p_1(n), k_j)$ and $\deg q_{1, j}< \deg p_1$.
	In this case $P'= \{q_{i,j}(n): (i,j)\in \{1,\cdots, d\}\times\{0,1,\cdots,l\}\} $. We still have $P'\subset \widetilde{SGP}$ and $\Phi(P')<\Phi(P)$. 
	
	Since $|k_j|>|k_{j-1}|+r(|k_{j-1}|)$ for $j=0, \cdots, l$,
	by Lemma \ref{lem-multi-q(ij)},
	$|A(q_{i, j})|>>1$ and $|A(q_{i,j}-q_{i',j'})|>>1$.

	By the inductive assumption, for $V_1^{(l)}, \cdots, V_d^{(l)}$.
	We have
	$$E=\{n\in \mathbb{Z}: V_1^{(l)}\cap \bigcap_{j=0}^{l}(T^{-q_{1, j}( n)}V_1^{(l)} \cap \cdots \cap T^{-q_{d, j}( n)}V_d^{(l)}) \neq \emptyset\} \cap (\tilde{C}\cap C_1)$$
	is syndetic.
	
	Let $m\in E$. We have $p_{i, j}(m) =q_{i, j}(m)$ since $m\in C_1$.
	Then there is some $x_m\in V_1^{(l)}$ such that
	$$T^{p_{i, j}( m)}x_m \in V_i^{(l)} \
	\text{for all} \ 1\leq i \leq d \ \text{and} \ 0\leq j\leq l.$$
	Clearly, there is some $y_m \in X$ such that $y_m=T^{-p_1(m)}x_m$.
	Since $X=\cup_{i=0}^{l}T^iU$, there is some $b_m\in \{0, 1, \cdots, l\}$
	such that
	$T^{b_m}z_m=y_m$ for some $z_m\in U$.
	Thus for each $i=1, \cdots, d$
	\begin{align*}
	T^{p_i(m+k_{b_m})}z_m
	&=T^{p_i(m+k_{b_m})}T^{-b_m}T^{-p_1(m)}x_m \\
	&=T^{p_i(k_{b_m})}T^{-b_m}T^{p_i(m+k_{b_m})}T^{-p_i(k_{b_m})}T^{-p_1(m)}x_m \\
	&=T^{p_i(k_{b_m})}T^{-b_m}T^{p_{i,{b_m}}(m)}x_m\\
	&\in T^{p_i(k_{b_m})}T^{-b_m} V_i^{(l)} \subset V_i.
	\end{align*}
	That is,
	$$z_m\in U\cap T^{-p_1(n)}V_1 \cap  \cdots \cap T^{-p_d(n)}V_d,$$
	where $n=m+k_{b_m}\in N$.
	
	Note that
	$k_{b_m}\in \widetilde{C}$
	implies
	$$  \{g_j(k_{b_m})\} \in (- \frac{\varepsilon}{2}, \frac{\varepsilon}{2} ),$$
	and
	$m\in C_1''$ implies
	$$  \{g_j(m+k_{b_m})\} \in ( \{g_j(k_{b_m})\} -\frac{\varepsilon}{2}, \{g_j(k_{b_m})\} +\frac{\varepsilon}{2} ) \subset (-\varepsilon, \varepsilon),$$
	for all $j=1, \cdots, t$.
	That is $m+k_{b_m}\in C$.
	
	Thus
	$$N \cap C\supset \{m+k_{b_m}: m\in E\}$$
	is a syndetic set.
	
	For every $P$, the induction will stop after finitely many steps, 	by induction the proof is completed.

\end{proof}

\medskip
\begin{proof}[Proof of Theorem \ref{thm general}.]
	By Lemma \ref{lem-equivalent}, it suffices to prove the moreover part of the theorem.
	Let $p_1, \cdots, p_d\in \mathcal{G}$. Then by Lemma \ref{lem sgp_gp}, there exists $h_i(n) \in \widetilde{SGP}$, $i=1,2,\dots,d$ and  $C_1=C(\delta,q_1,\cdots,q_k)$ such that $$p_i(n)=_{C_1}h_i(n), i=1,2,\dots,d.$$	
	Set $$N_1=\{n \in \mathbb{N}: U\cap T^{-h_1(n)}V_1\cap  \cdots \cap T^{-h_d(n)}V_d \neq \emptyset \},$$
	by Theorem \ref{thm general sgp}, $N_1\cap (C\cap C_1)$ is syndetic. Since for any $n\in N_1\cap (C\cap C_1) \subset C_1$, $p_i(n)=h_i(n),i=1,2,\cdots,d$, 
	$$n \in N=\{n \in \mathbb{N} : U \cap T^{-p_1(n)}V_1\cap \cdots T^{-p_d(n)}V_d \neq \emptyset \}.$$
	This implies	
	$$N_1 \cap (C \cap C_1) \subset N \cap C,$$
	hence $N\cap C$ is syndetic.
\end{proof}

\end{document}